\documentstyle[11pt]{article}
\begin{document}

\title{Integrable systems in projective differential geometry}
\author{{\Large Ferapontov E.V. } \\
    Department of Mathematical Sciences \\
    Loughborough University \\
    Loughborough, Leicestershire LE11 3TU \\
    United Kingdom \\
    and \\
    Centre for Nonlinear Studies \\
    Landau Institute of Theoretical Physics\\
    Academy of Science of Russia, Kosygina 2\\
    117940 Moscow, GSP-1, Russia\\
    e-mail: {\tt fer@landau.ac.ru}
}
\date{}
\maketitle

\newtheorem{theorem}{Theorem}

\pagestyle{plain}

\maketitle

\begin{abstract}

Some of the most important classes of surfaces in projective 3-space are
reviewed: these are isothermally asymptotic surfaces, projectively applicable
surfaces, surfaces of Jonas,
projectively minimal surfaces, etc. It is demonstrated that the
corresponding projective "Gauss-Codazzi" equations reduce to integrable systems
which are quite familiar from the modern soliton theory and coincide with the
stationary flows in the Davey-Stewartson and Kadomtsev-Petviashvili
hierarchies,
equations of the Toda lattice, etc. The corresponding Lax pairs
can be obtained by inserting a spectral parameter in the equations of the
Wilczynski moving frame.

\bigskip


\bigskip

\end{abstract}

\newpage
\section{Introduction}

Projective differential geometry of surfaces $M^2$ in $P^3$ has been
extensively
developed in the first half of the century in the works of Wilczynski, Fubini,
\^Cech, Cartan,  Tzitzeica, Demoulin, Rozet, Godeaux, Lane, Eisenhart, Finikov,
Bol and many others.
These investigations culminated in a number of beautifull geometric
constructions
and a list of particularly interesting classes of surfaces which now seem
to be quite deeply
forgotten. Our aim here is to give a brief review of surfaces in projective
3-space with the emphasize on the nonlinear equations underlying them. In
particular, we  demonstrate that
\medskip \par
\noindent {\bf isothermally asymptotic surfaces} are described by the
stationary modified
Veselov-Novikov (mVN) equation (sect.3);
\medskip \par
\noindent {\bf projectively applicable surfaces} (which naturally fall into the
classes
of the so-called "surfaces $R_0$" and "surfaces $R$") correspond to the
stationary limits of the second flow in the Kadomtsev-Petviashvili (KP)
hierarchy and the stationary Davey-Stewartson (DS) system, respectively
(sect.4);
\medskip \par
\noindent {\bf surfaces of Jonas} are related to a stationary limit of the
fourth order
flow in the DS hierarchy (sect.5);
\medskip \par
\noindent {\bf projectively minimal surfaces} are governed by an integrable
system which, in a certain
limit, reduces to the coupled Tzitzeica system, being a reduction of
periodic
Toda lattice of period 6 (sect.6);
\medskip \par
\noindent {\bf Surfaces with asymptotic lines in linear complexes}
correspond to a linearizable system (sect.7).
\medskip

Thus, all of the most important classes of surfaces in projective
differential geometry are described by integrable systems which are an
object of the modern soliton theory. This is by no means surprising: the
reason for differential geometers to introduce a particular class of surfaces
has always been the existence of certain "nice geometric transformations"
allowing the construction of  new surfaces from the given ones.
On the other hand, the existence of B\"acklund transformations is a
standard hint on the integrability. Thus, in a sense, classical differential
geometry (here we mean geometry of surfaces in 3-space) can be viewed as a
chapter
of the theory of integrable systems. Most of the
B\"acklund transformations derived in the  "solitonic" context can be found
in the old geometric textbooks (probably, stated in geometric terms, without
an explicit coordinate formulae). In sect.7 we review the
construction of congruences W, which are the main source of
B\"acklund transformations.
It requires a solution of certain Dirac equation on the surface $M^2$.
Since Dirac operator plays a role of the Lax operator of the (2+1)-dimensional
DS hierarchy, this construction clarifies the relationship between surfaces in
projective differential geometry and stationary flows of the DS hierarchy.

In our approach we make use of the Wilczynski moving frame which proves to be
the most convenient tool for studying surfaces in projective 3-space. For
instance, the Lax pairs for  systems governing particular classes of surfaces
can be obtained by inserting an appropriate
spectral parameter in the equations of motion of the Wilczynski frame.

\section{Surfaces in projective differential geometry}

Based on \cite{Wilczynski} (see also \cite{Fer8}, \cite{Sasaki1},
\cite{Sasaki2}), let us briefly recall
the standard way of defining surfaces $M^2$
in projective space $P^3$ in terms of solutions of a linear system
\begin{equation}
\begin{array}{c}
{\bf r}_{xx}=\beta \ {\bf r}_y+\frac{1}{2}(V-\beta_y) \ {\bf r} \\
\ \\
{\bf r}_{yy}=\gamma \ {\bf r}_x+\frac{1}{2}(W-\gamma_x) \ {\bf r}
\end{array}
\label{r}
\end{equation}
where $\beta, \gamma, V, W$ are functions of $x$ and  $y$. If we
cross-differentiate (\ref{r}) and assume
${\bf r}, {\bf r}_x, {\bf r}_y, {\bf r}_{xy}$ to be independent, we arrive at
the compatibility
conditions \cite[p.\ 120]{Lane}
\begin{equation}
\begin{array}{c}
\beta_{yyy}-2\beta_yW-\beta W_y=
\gamma_{xxx}-2\gamma_xV-\gamma V_x\\
\ \\
W_x=2\gamma \beta_y+\beta \gamma_y \\
\ \\
V_y=2\beta \gamma_x+\gamma \beta_x.
\end{array}
\label{GC1}
\end{equation}
For any fixed
$\beta, \gamma, V, W$ satisfying (\ref{GC1}) the linear system
(\ref{r})  is compatible
and possesses a solution
${\bf r}=(r^0, r^1, r^2, r^3)$ where $r^i(x, y)$
 can be regarded as homogeneous coordinates of a
surface in projective space $P^3$. One may think of
$M^2$ as a surface in a three-dimensional  space with position vector
${\bf R}=(r^1/r^0, r^2/r^0, r^3/r^0)$.
If we choose any other solution
$\tilde{{\bf r}}=(\tilde r^0, \tilde r^1, \tilde r^2, \tilde r^3)$
of the same system
(\ref{r}) then the corresponding surface
$\tilde M^2$ with position vector
$\tilde{{\bf R}}=(\tilde r^1/\tilde r^0, \tilde r^2/\tilde r^0,
\tilde r^3/\tilde r^0)$
constitutes a projective transform of $M^2$ so that any fixed
$\beta, \gamma, V, W$ satisfying (\ref{GC1}) define a surface $M^2$ uniquely up
to
projective equivalence. Moreover, a simple calculation yields
\begin{equation}
\begin{array}{c}
{\bf R}_{xx}=\beta \ {\bf R}_y+ a \ {\bf R}_x \\
{\bf R}_{yy}=\gamma \ {\bf R}_x+ b \ {\bf R}_y
\end{array}
\label{affine}
\end{equation}
$(a=-2 r^0_x/r^0,\,b=-2 r^0_y/r^0)$ which implies that $x, y$ are asymptotic
coordinates of the surface $M^2$. System (\ref{affine}) can be viewed
as an "affine gauge" of system (\ref{r}).
In what follows, we assume that our surfaces are
hyperbolic and the corresponding asymptotic coordinates $x, y$ are
real.\footnote{The elliptic case is dealt with in an analogous manner by
regarding $x,y$ as complex conjugates.}
Since equations~(\ref{GC1})
specify a surface uniquely up to projective equivalence, they can be viewed as
the `Gauss-Codazzi' equations in projective geometry.

Even though the coefficients $\beta, \gamma, V, W$ define a surface $M^2$
uniquely up to projective equivalence via (\ref{r}),
it is not entirely correct
to regard $\beta, \gamma, V, W$ as projective invariants. Indeed, the
asymptotic
coordinates $x, y$ are only defined up to an arbitrary reparametrization
of the form
\begin{equation}
x^*=f(x), ~~~~ y^*=g(y)
\label{newxy}
\end{equation}
which induces a scaling of the surface vector according to
\begin{equation}
{\bf r}^*=\sqrt {f'(x)g'(y)}~{\bf r}.
\label{newr}
\end{equation}
Thus \cite[p.\ 1]{Bol1}, the form of
equations (\ref{r}) is preserved by the above transformation
with the new coefficients $\beta^*, \gamma^*, V^*, W^*$ given by
\begin{equation}
\begin{array}{c}
\beta^{*}=\beta g'/(f')^2, ~~~~  V^{*}(f')^2=V+S(f)\\
\ \\
\gamma^{*}=\gamma f'/(g')^2, ~~~~ W^{*}(g')^2=W+S(g),
\end{array}
\label{new}
\end{equation}
where $S(\,\cdot\,)$ is the  Schwarzian derivative, that is
$$
S(f)=\frac{f'''}{f'} - \frac{3}{2} \left(\frac{f''}{f'}\right)^2.
$$
The transformation formulae (\ref{new}) imply that the symmetric 2-form
$$
2 \beta \gamma\,dxdy
$$
and the conformal class of the cubic form
$$
\beta \,dx^3+\gamma \,dy^3
$$
are absolute projective invariants. They are known as the projective metric and
the Darboux cubic form, respectively, and play an important role in projective
differential geometry. In particular, they define a `generic' surface
uniquely up to projective equivalence. The vanishing of the Darboux cubic form
is characteristic for quadrics: indeed, in this case $\beta = \gamma =0$ so that
asymptotic curves of both families are straight lines. The vanishing of the
projective metric (which is equivalent to either $\beta =0$ or $\gamma =0$)
characterises ruled surfaces. In what follows we exclude these two degenerate
situations and require $\beta \ne 0$, $\gamma \ne 0$.

One can also define projectively invariant differentials
$$
\omega^1=(\gamma \beta^2)^{\frac{1}{3}}\, dx, ~~~~
\omega^2=(\beta \gamma^2)^{\frac{1}{3}}\, dy, ~~~~
$$
so that $\beta \gamma \,  dxdy=\omega^1 \omega^2$. With the help of
$\omega^1, ~ \omega^2$ one can define projectively invariant differentiation.
However, we will not take advantage of it in what follows.

Using (\ref{newxy})-(\ref{new}), one can  verify that the four
points
\begin{equation}
\begin{array}{c}
{\bf r}, ~~~
{\bf r}_1={\bf r}_x-\frac{1}{2}\frac{\gamma_x}{\gamma}{\bf r}, ~~~
{\bf r}_2={\bf r}_y-\frac{1}{2}\frac{\beta_y}{\beta}{\bf r}, \\
\ \\
{\mbox{\boldmath $\eta$}}={\bf r}_{xy}-\frac{1}{2}\frac{\gamma_x}{\gamma}{\bf
r}_y-
\frac{1}{2}\frac{\beta_y}{\beta}{\bf r}_x+
\left(\frac{1}{4}\frac{\beta_y\gamma_x}{\beta \gamma} -
\frac{1}{2}{\beta \gamma}\right){\bf r} \\
\end{array}
\label{frame}
\end{equation}
are defined in an invariant way, that is under the transformation formulae
\mbox{(\ref{newxy})-(\ref{new})}
they acquire a nonzero multiple which does not
change them as points in projective space $P^3$. These points form the
vertices of the so-called Wilczynski moving tetrahedral \cite{Bol1},
\cite{Finikov37}, \cite{Wilczynski}.
Since the lines passing through ${\bf r}, {\bf r}_1$ and ${\bf r}, {\bf r}_2$
are
tangential to the $x$- and $y$-asymptotic curves,
respectively, the three points ${\bf r}, {\bf r}_1, {\bf r}_2$
span the tangent plane of the surface $M^2$.
The line through ${\bf r}_1, {\bf r}_2$ lying in the tangent
plane is
known as the directrix of Wilczynski of the second kind. The line through
${\bf r}, {\mbox{\boldmath $\eta$}}$ is transversal to $M^2$ and is known as
the directrix of
Wilczynski of the first kind. It plays the role of a projective `normal'.
We stress that in projective differential geometry there exists no unique
choice of an invariant normal. This is in contrast with Euclidean
and affine geometries in which the normal is canonically defined.
Some of the best-known and most-investigated normals are those of
Wilczynski, Fubini, Green, Darboux, Bompiani and Sullivan \cite[p.\ 35]{Bol1}
with the directrix of Wilczynski being the most commonly used. It is known that
the normal of Wilczynski intersects the tangent Lie quadric of the surface
$M^2$ at exactly two points ${\bf r}$ and ${\mbox{\boldmath $\eta$}}$ so that
both points lie on the
Lie quadric and are canonically defined. The Wilczynski tetrahedral proves to
be the most convenient tool in projective differential geometry.

Using (\ref{r}) and~(\ref{frame}),
we easily derive for ${\bf r}, {\bf r}_1, {\bf r}_2, {\mbox{\boldmath $\eta$}}$
the linear equations \cite[p.\ 42]{Finikov37}
\begin{equation}
\begin{array}{c}
\left(\begin{array}{c}
{\bf r}\\
{\bf r}_1\\
{\bf r}_2\\
{\mbox{\boldmath $\eta$}}
\end{array}\right)_x=
\left(\begin{array}{cccc}
\frac{1}{2}\frac{\gamma_x}{\gamma}& 1&0&0\\
\frac{1}{2}b & -\frac{1}{2}\frac{\gamma_x}{\gamma}& \beta&0\\
\frac{1}{2}k&0&\frac{1}{2}\frac{\gamma_x}{\gamma}&1\\
\frac{1}{2}\beta
a&\frac{1}{2}k&\frac{1}{2}b&-\frac{1}{2}\frac{\gamma_x}{\gamma}
\end{array}\right)
\left(\begin{array}{c}
{\bf r}\\
{\bf r}_1\\
{\bf r}_2\\
{\mbox{\boldmath $\eta$}}
\end{array}\right)\\
\ \\
\left(\begin{array}{c}
{\bf r}\\
{\bf r}_1\\
{\bf r}_2\\
{\mbox{\boldmath $\eta$}}
\end{array}\right)_y=
\left(\begin{array}{cccc}
\frac{1}{2}\frac{\beta_y}{\beta}& 0&1&0\\
\frac{1}{2}l & \frac{1}{2}\frac{\beta_y}{\beta}& 0&1\\
\frac{1}{2}a&\gamma &-\frac{1}{2}\frac{\beta_y}{\beta}&0\\
\frac{1}{2}\gamma b&\frac{1}{2}a
&\frac{1}{2}l&-\frac{1}{2}\frac{\beta_y}{\beta}
\end{array}\right)\left(
\begin{array}{c}
{\bf r}\\
{\bf r}_1\\
{\bf r}_2\\
{\mbox{\boldmath $\eta$}}
\end{array}\right),
\end{array}
\label{Wilczynski}
\end{equation}
where we introduced the notation
\begin{equation}
\begin{array}{c}
k=\beta \gamma - (\ln \beta)_{xy}, ~~~~ l=\beta \gamma - (\ln \gamma)_{xy},\\
\ \\
a=W-(\ln \beta)_{yy}-\frac{1}{2}(\ln \beta)_y^2, ~~~~
b=V-(\ln \gamma)_{xx}-\frac{1}{2}(\ln \gamma)_x^2.
\end{array}
\label{klab}
\end{equation}
Under the transformations (\ref{newxy})-(\ref{new})
these quantities transform as follows
\begin{equation}
\begin{array}{c}
k^{*}= k/f'g', ~~~~ l^{*}= l/f'g',\\
\ \\
a^{*}=a/(g')^2, ~~~~ b^{*}=b/(f')^2,
\end{array}
\label{(newklab)}
\end{equation}
and give rise to the projectively invariant quadratic form
$$
b\, dx^2+a\, dy^2
$$
and the quartic form
$$
a\beta^2dx^4+b\gamma^2dy^4.
$$

The compatibility conditions of equations (\ref{Wilczynski}) imply
\begin{equation}
\begin{array}{c}
(\ln \beta)_{xy}=\beta \gamma -k, ~~~~ (\ln \gamma)_{xy}=\beta \gamma -l,\\
\ \\
a_x=k_y+\frac{\beta_y}{\beta}k, ~~~~ b_y=l_x+\frac{\gamma_x}{\gamma}l,\\
\ \\
\beta a_y+2a\beta_y=\gamma b_x+2b\gamma_x,
\end{array}
\label{GC2}
\end{equation}
which is just the equivalent form of the projective "Gauss-Codazzi" equations
(\ref{GC1}).

Equations (\ref{Wilczynski}) can be rewritten in the Pl\"ucker coordinates.
For a convenience of the reader we briefly recall this construction.
Let us consider a line $l$ in $P^3$ passing through the points ${\bf a}$ and
${\bf b}$  with the
homogeneous coordinates ${\bf a}=(a^0:a^1:a^2:a^3)$ and ${\bf
b}=(b^0:b^1:b^2:b^3)$. With
the line $l$ we associate a point ${\bf a}\wedge {\bf b}$ in projective space
$P^5$ with
the homogeneous coordinates
$$
{\bf a}\wedge {\bf b}=(p_{01}:p_{02}:p_{03}:p_{23}:p_{31}:p_{12}),
$$
where
$$
p_{ij}=\det
\left(\begin{array}{cc}
a^i & a^j \\
b^i & b^j
\end{array}\right).
$$
The coordinates $p_{ij}$ satisfy the well-known quadratic Pl\"ucker relation
\begin{equation}
p_{01}\, p_{23}+p_{02}\, p_{31}+p_{03}\, p_{12}=0.
\label{quadric}
\end{equation}
Instead of ${\bf a}$ and ${\bf b}$ we may consider an arbitrary linear
combinations thereof
without changing ${\bf a}\wedge{\bf b}$ as a point in $P^5$.
Hence, we arrive at the well-defined
Pl\"ucker coorrespondence $l({\bf a},{\bf b})\to {\bf a}\wedge {\bf b}$
between lines in
$P^3$ and points on the Pl\"ucker quadric in $P^5$.
Pl\"ucker correspondence plays an important role in the projective
differential geometry of surfaces and often sheds some new light on those
properties of surfaces which are not `visible' in $P^3$ but acquire a precise
geometric meaning only in $P^5$. Thus, let us consider a surface
$M^2\in P^3$ with the Wilczynski tetrahedral
${\bf r}, {\bf r}_1, {\bf r}_2, {\mbox{\boldmath $\eta$}} $
satisfying equations (\ref{Wilczynski}).
Since the two pairs of points ${\bf r}, {\bf r}_1$ and ${\bf r}, {\bf r}_2$
generate two lines in $P^3$ which are tangential
to the $x$- and $y$-asymptotic curves, respectively, the formulae
$$
{\cal U}={\bf r} \wedge {\bf r}_1,\quad {\cal V}={\bf r} \wedge {\bf r}_2
$$
define the images of these lines under the Pl\"ucker embedding. Hence, with any
surface $M^2\in P^3$ there are canonically associated two surfaces
${\cal U}(x, y)$ and ${\cal V}(x, y)$ in $P^5$ lying on the
Pl\"ucker quadric (\ref{quadric}).
In view of the formulae
$$
{\cal U}_x=\beta \, {\cal V},\quad {\cal V}_y=\gamma \, {\cal U},
$$
we conclude that the line in $P^5$ passing through a pair of points
$({\cal U}, {\cal V})$ can also be generated by the pair of points
$({\cal U}, {\cal U}_x)$ (and hence is
tangential to the $x$-coordinate line on the surface ${\cal U}$) or by a pair
of
points $({\cal V}, {\cal V}_y)$ (and hence is tangential to the $y$-coordinate
line on the
surface ${\cal V}$).
Consequently, the surfaces ${\cal U}$ and ${\cal V}$ are two focal surfaces of
the congruence
of straight
lines $({\cal U}, {\cal V})$ or, equivalently,
${\cal V}$ is the Laplace transform of ${\cal U}$ with
respect to $x$ and ${\cal U}$ is the
Laplace transform of ${\cal V}$ with respect to $y$.
We emphasize that the $x$- and
$y$-coordinate lines on the surfaces ${\cal U}$ and ${\cal V}$
are not asymptotic but
conjugate. Continuation of the Laplace sequence in both directions, that is
taking the $x$-transform of ${\cal V}$, the $y$-transform of ${\cal U}$, etc.,
leads, in the
generic case, to an infinite Laplace sequence in $P^5$ known as the Godeaux
sequence of a surface $M^2$ \cite[p.\ 344]{Bol1}. The surfaces of the Godeaux
sequence carry important geometric information about the surface $M^2$
itself.

The case of a closed, i.e. periodic Godeaux sequence is particularly
interesting.
It turns out, that the only surfaces $M^2\in P^3$ for which the Godeaux
sequence
is of period 6 (the value 6 turns out to be the least possible) are the
surfaces of
Demoulin \cite[p.\ 360]{Bol1} -- see sect.4.

Introducing
$$
\begin{array}{c}
{\cal A} ={\bf r}_2\wedge {\bf r}_1+{\bf r}\wedge {\mbox{\boldmath $\eta$}},
~~~
{\cal B} ={\bf r}_1\wedge {\bf r}_2+{\bf r}\wedge {\mbox{\boldmath $\eta$}}, \\
\ \\
{\cal P} = 2\, {\bf r}_2\wedge {\mbox{\boldmath $\eta$}}, ~~~
{\cal Q} = 2\, {\bf r}_1\wedge {\mbox{\boldmath $\eta$}},
\end{array}
$$
we arrive at the following equations for the Pl\"ucker coordinates:

\begin{equation}
\begin{array}{c}
\left(\begin{array}{c}
{\cal U}\\
{\cal A}\\
{\cal P}\\
{\cal V}\\
{\cal B}\\
{\cal Q}
\end{array}\right)_x=
\left(\begin{array}{cccccc}
0 & 0 & 0 & \beta & 0 & 0\\
k & 0 & 0 & 0 & 0 & 0\\
0 & k & 0 & -\beta a & 0 & 0\\
0 & 0 & 0 & \frac{\gamma_x}{\gamma} & 1 & 0\\
0 & 0 & 0 & b & 0 & 1\\
-\beta a & 0 & \beta & 0 & b &-\frac{\gamma_x}{\gamma}
\end{array}\right)
\left(\begin{array}{c}
{\cal U}\\
{\cal A}\\
{\cal P}\\
{\cal V}\\
{\cal B}\\
{\cal Q}
\end{array}\right)\\
\ \\
\left(\begin{array}{c}
{\cal U}\\
{\cal A}\\
{\cal P}\\
{\cal V}\\
{\cal B}\\
{\cal Q}
\end{array}\right)_y=
\left(\begin{array}{cccccc}
\frac{\beta_y}{\beta} & 1 & 0 & 0 & 0 & 0\\
a & 0 & 1 & 0 & 0 & 0\\
0 & a & -\frac{\beta_y}{\beta} & -\gamma b & 0 & \gamma\\
\gamma & 0 & 0 & 0 & 0 & 0\\
0 & 0 & 0 & l & 0 & 0\\
-\gamma b & 0 & 0 & 0 & l & 0
\end{array}\right)
\left(\begin{array}{c}
{\cal U}\\
{\cal A}\\
{\cal P}\\
{\cal V}\\
{\cal B}\\
{\cal Q}
\end{array}\right)
\end{array}
\label{UAPVBQ}
\end{equation}
Equations (\ref{UAPVBQ}) are consistent with the following table
of scalar products:
\begin{equation}
({\cal U}, {\cal P})=-1, ~~~ ({\cal A}, {\cal A})=1, ~~~
({\cal V}, {\cal Q})=1, ~~~ ({\cal B}, {\cal B})=-1,
\label{table}
\end{equation}
all other scalar products being equal to zero. This defines a scalar
product of the signature (3, 3) which is the same as that of the quadratic form
(\ref{quadric}).

Different types of surfaces can be defined by imposing additional constraints
on $\beta$, $\gamma$, $V$, $W$ (respectively, $\beta, \gamma, k, l, a, b$),
so that, in a sense, projective differential geometry is the
theory of (integrable) reductions of the underdetermined system (\ref{GC1})
(respectively, (\ref{GC2})).
Although the three linear systems (\ref{r}),
(\ref{Wilczynski}) and (\ref{UAPVBQ}) are in fact equivalent, some of them
prove
to be more suitable for studying particular classes of projective surfaces.

{\bf Remark.} Since the tangent plane of the surface $M^2$
is spanned by three points ${\bf r}, {\bf r }_1, {\bf r}_2$, the vector
${\bf r}^d= {\bf r}\wedge {\bf r}_1\wedge {\bf r}_2$ can be viewed as
a radius-vector of the dual surface.
A simple calculation yields
$$
\begin{array}{c}
{\bf r}^d_{xx}=-\beta \ {\bf r}^d_y+\frac{1}{2}(V+\beta_y) \
{\bf r}^d \\
\ \\
{\bf r}^d_{yy}=-\gamma \ {\bf r}^d_x+\frac{1}{2}(W+\gamma_x) \
{\bf r}^d
\end{array}
$$
implying that the passage to the dual surface is equivalent to a simple
change of signs:
$\beta, \gamma, V, W \to -\beta, -\gamma, V, W$.
This transformation is obviously a discrete symmetry of
equations (\ref{GC1}).

\section {Isothermally asymptotic surfaces}

These surfaces are specified by the condition
$$
\left(\ln {\frac{\beta}{\gamma}}\right)_{xy}=0,
$$
which, in view of the transformation formulae (\ref{new}),
reduces to $\beta=\gamma$ after a suitable choice
of coordinates $x, y$. In this case
equations (\ref{GC1}) assume the form
of the stationary modified Veselov-Novikov (mVN) equation
\begin{equation}
\begin{array}{c}
\beta_{yyy}-2\beta_yW-\beta W_y =
\beta_{xxx}-2\beta_xV-\beta V_x, \\
\, \\
W_x=\frac{3}{2}(\beta^2)_y \\
\ \\
V_y=\frac{3}{2}(\beta^2)_x.
\end{array}
\label{smVN}
\end{equation}
In this case one can introduce a spectral parameter $\lambda$ in equations
(\ref{UAPVBQ}) without violating their compatibility:
\begin{equation}
\begin{array}{c}
\left(\begin{array}{c}
{\cal U}\\
{\cal A}\\
{\cal P}\\
{\cal V}\\
{\cal B}\\
{\cal Q}
\end{array}\right)_x=
\left(\begin{array}{cccccc}
0 & 0 & 0 & \beta & 0 & 0\\
k & 0 & 0 & 0 & 0 & 0\\
0 & k & 0 & -\beta a & 0 & 0\\
0 & 0 & 0 & \frac{\beta_x}{\beta} & 1 & 0\\
0 & 0 & 0 & b & 0 & 1\\
-\beta a & 0 & \beta & \lambda & b &-\frac{\beta_x}{\beta}
\end{array}\right)
\left(\begin{array}{c}
{\cal U}\\
{\cal A}\\
{\cal P}\\
{\cal V}\\
{\cal B}\\
{\cal Q}
\end{array}\right)\\
\ \\
\left(\begin{array}{c}
{\cal U}\\
{\cal A}\\
{\cal P}\\
{\cal V}\\
{\cal B}\\
{\cal Q}
\end{array}\right)_y=
\left(\begin{array}{cccccc}
\frac{\beta_y}{\beta} & 1 & 0 & 0 & 0 & 0\\
a & 0 & 1 & 0 & 0 & 0\\
-\lambda & a & -\frac{\beta_y}{\beta} & -\beta b & 0 & \beta\\
\beta & 0 & 0 & 0 & 0 & 0\\
0 & 0 & 0 & l & 0 & 0\\
-\beta b & 0 & 0 & 0 & l & 0
\end{array}\right)
\left(\begin{array}{c}
{\cal U}\\
{\cal A}\\
{\cal P}\\
{\cal V}\\
{\cal B}\\
{\cal Q}
\end{array}\right)
\end{array}
\label{UAPVBQ1}
\end{equation}
Rewriting (\ref{UAPVBQ1}) in terms of ${\cal U}, {\cal V}$ (that is, expressing
${\cal A}, {\cal B}, {\cal P}, {\cal Q}$ through
${\cal U}, {\cal U} _y, {\cal U} _{yy}$ and ${\cal V}, {\cal V} _x,
{\cal V} _{xx}$), we arrive at
\begin{equation}
\begin{array}{c}
{\cal U}_x=\beta \ {\cal V} \\
{\cal V}_y=\beta \ {\cal U} \\
\\
\lambda {\cal U}={\cal U}_{xxx}- {\cal U}_{yyy}+2 \ W \ {\cal U}_y-
3\ \beta_x \ {\cal V}_x+ W_y \ {\cal U} -2\ \beta \ V \ {\cal V} \\
\\
\lambda {\cal V}= {\cal V}_{xxx}- {\cal V}_{yyy}-2 \ V \ {\cal V}_x +
3 \ \beta _y \ {\cal U}_y -  V_x \ {\cal V} +2 \ \beta \   W \ {\cal U}. \\
\end{array}
\label{slmVN}
\end{equation}
which coincides with the stationary limit of the mVN linear problem
(the compatibility conditions of  (\ref{slmVN}) coincide with
(\ref{smVN})).

We recall that the (2+1)-dimensional mVN equation
\begin{equation}
\begin{array}{c}
\beta_t= \beta_{xxx}-\beta_{yyy}-2 \ \beta_x \ V+2 \ \beta_y \ W-\beta \ V_x+
\beta \ W_y \\
\\
W_x=\frac{3}{2}(\beta^2)_y \\
\\
V_y=\frac{3}{2}(\beta^2)_x \\
\end{array}
\label{mVN}
\end{equation}
was introduced in \cite{Bogdanov} and is associated with the two-dimensional
Dirac operator
\begin{equation}
\begin{array}{c}
{\cal U}_x=\beta \ {\cal V} \\
{\cal V}_y=\beta \ {\cal U} \\
\\
{\cal U}_t={\cal U}_{xxx}- {\cal U}_{yyy}+2 \ W \ {\cal U}_y-
3\ \beta_x \ {\cal V}_x+ W_y \ {\cal U} -2\ \beta \ V \ {\cal V} \\
\\
{\cal V}_t= {\cal V}_{xxx}- {\cal V}_{yyy}-2 \ V \ {\cal V}_x +
3 \ \beta _y \ {\cal U}_y -  V_x \ {\cal V} +2 \ \beta \   W \ {\cal U}. \\
\end{array}
\label{lmVN}
\end{equation}
Linear system (\ref{slmVN}) results after a  substitution
${\cal U}_t \to \lambda {\cal U}$, ~ ${\cal V}_t \to \lambda {\cal V}$
which is a standard way to introduce a spectral parameter
in the stationary problem.

Isothermally asymptotic surfaces are attributed to Fubini
\cite{Fubini}  and
can be equivalently defined by any of the following geometric properties:

- The 3-web, formed by asymptotic curves and Darboux's curves, is
hexagonal. (Darboux's curves are the zero curves of the
Darboux cubic form $\beta \,  dx^3+\gamma \, dy^3$).

- Isothermally asymptotic surfaces are the focal surfaces of  special
$W$-congruences, preserving Darboux's curves (in fact, this
property implies a B\"acklund transformation for isothermally asymptotic
surfaces -- see sect.7).

Examples of isothermally asymptotic surfaces include arbitrary quadrics and
cubics, quartics of Kummer, projective transforms of affine spheres and
rotation surfaces:

\bigskip

{\bf Quadrics} correspond to the trivial solution $\beta=0, ~ W=W(y), ~
V=V(x)$.
\bigskip

{\bf Projective transforms of rotation surfaces} are specified by
$\beta=\beta(x+y), ~ W=V=\frac{3}{2}\beta^2+c$ where $\beta $
is an arbitrary function of $(x+y)$ and
$c$ is an arbitrary constant. For $c>0$ these are indeed projective transforms
 of surfaces $z=f(x^2+y^2)$, while the cases $c=0$ and $c<0$
correspond to projective transforms of
 surfaces $z=f(x^2+y)$ and $z=f(x^2-y^2)$, respectively.
Travelling-wave solutions $\beta (x+cy)$
of equation (\ref{smVN}) correspond to surfaces, which are invariant
under one-parameter groups of projective transformations. In the case $c\ne 1$
the function $\beta$ is no longer arbitrary and can be expressed in elliptic
functions
(compare with \cite{Lingenberg1}).

\bigskip

{\bf Cubic surfaces} are specified by the following additional constraints in
(\ref{smVN}):
\begin{equation}
V=-\frac{1}{2}(\ln \beta)_{xx}+\frac{1}{8}(\ln \beta)_x^2+\frac{5}{2}\beta_y,
~~~
W=-\frac{1}{2}(\ln \beta)_{yy}+\frac{1}{8}(\ln \beta)_y^2+\frac{5}{2}\beta_x
\label{1}
\end{equation}
\cite{Lane1}, see also \cite{Lane}, p.131. With these $V,W$ equations
(\ref{smVN}) imply
$$
\left(\frac{(\ln \beta)_{xy}}{\sqrt \beta}+4\beta\sqrt \beta\right)_y=
5\frac{\beta_{xx}}{\sqrt \beta},
~~~
\left(\frac{(\ln \beta)_{xy}}{\sqrt \beta}+4\beta\sqrt \beta\right)_x=
5\frac{\beta_{yy}}{\sqrt \beta}.
$$
Integration of these equations for $\beta$ would provide a 4-parameter family
of exact solutions of equation (\ref{smVN}): indeed, up to projective
equivalence cubics in $P^3$ depend on 4 essential parameters.
\bigskip

{\bf The Roman surface of Steiner} is a rational quartic in $P^3$ with the
equation
$$
(x^2+y^2+z^2-1)^2=((z-1)^2-2x^2)((z+1)^2-2y^2)
$$
owing it's name to   Steiner who investigated this surface
in Rome in 1844. Besides quadrics and ruled cubic surfaces
the Roman surface of Steiner is the only surface in $P^3$
possessing infinitely many conic sections
through any of it's points. This result was announced
several times: by Moutard in 1865,  Darboux in 1880 and Wilczynski in 1908 (see
\cite{Wilczynski}, 1909 for historical remarks). The interest to the
Roman surface of Steiner in projective differential geometry is due to
the remarkable
construction of Darboux, relating with an arbitrary surface $M^2$ in $P^3$
and an arbitrary point $p$ on $M^2$ an osculating Roman surface of Steiner
which has the fourth order of tangency with $M^2$ at this point.
Analytically, the Roman surface of Steiner corresponds to the choice
\begin{equation}
V=-\frac{1}{2}(\ln \beta)_{xx}+\frac{1}{8}(\ln \beta)_x^2-\frac{5}{2}\beta_y,
~~~
W=-\frac{1}{2}(\ln \beta)_{yy}+\frac{1}{8}(\ln \beta)_y^2-\frac{5}{2}\beta_x,
\label{2}
\end{equation}
$$
(\ln \beta)_{xy}=\frac{4}{9}\beta^2
$$
(\cite{Bol1}, p.149-150)
implying upon substitution in (\ref{smVN})  the following equations for
$\beta$:
\begin{equation}
\begin{array}{c}
\beta_{xx}=-\frac{4}{3}\beta\beta_y\\
\\
\beta_{yy}=-\frac{4}{3}\beta\beta_x\\
\\
(\ln \beta)_{xy}=\frac{4}{9}\beta^2.
\end{array}
\label{Roman}
\end{equation}
These can be explicitely integrated:
$$
\beta^2=\frac{9}{4}\frac{f^{'}g^{'}}{(f+g)^2}
$$
where the functions $f(x)$ and $g(y)$ satisfy the ODE's
$$
(f^{'})^3=(a_0+a_1f+a_2f^2)^2, ~~~ (g^{'})^3=(a_0-a_1g+a_2g^2)^2.
$$
Here $a_i$ are arbitrary constants.
Under the transformation $(\beta, V, W)\to (-\beta, V, W)$ equations
(\ref{2}) transform to (\ref{1}). This means, that the dual of the Roman
surface of Steiner is a cubic, and hence the Roman surface itself is a quartic
of class 3 (\cite{Bol1}, p.150).

The Roman surface of Steiner belongs to
a broader class of isothermally asymptotic quartic surfaces known as

\bigskip

{\bf Quartics of Kummer} investigated by Kummer as singular surfaces of
quadratic line complexes. Analytically, the quartics of Kummer
are specified by the conditions
\begin{equation}
\begin{array}{c}
V=\frac{11}{8}(\ln \beta)_{xx}+2(\ln \beta)_x^2, ~~~
W=\frac{11}{8}(\ln \beta)_{yy}+2(\ln \beta)_y^2, \\
\ \\
(\ln \beta)_{xy}=\frac{4}{9}\beta^2
\end{array}
\label{3}
\end{equation}
(\cite{Bol1}, p.231). Substituting these $V, W$ in  (\ref{smVN}) we arrive at
\begin{equation}
\left(\frac{1}{\beta^2}(\beta^2(\beta^2)_y)_y\right)_y=
\left(\frac{1}{\beta^2}(\beta^2(\beta^2)_x)_x\right)_x.
\label{p}
\end{equation}
With
$$
\beta^2=\frac{9}{4}\frac{f^{'}g^{'}}{(f+g)^2}
$$
(which is a general solution of the Liouville equation $(\ref{3})_3$)
equations (\ref{p}) can be rewritten in the form
\begin{equation}
\begin{array}{c}
\frac{1}{3}(f+g)^3
\left(\frac{d^3((f^{'})^3)}{df^3}-\frac{d^3((g^{'})^3)}{dg^3}\right)-
4(f+g)^2
\left(\frac{d^2((f^{'})^3)}{df^2}-\frac{d^2((g^{'})^3)}{dg^2}\right)+ \\
\ \\
20(f+g)
\left(\frac{d((f^{'})^3)}{df}-\frac{d((g^{'})^3)}{dg}\right)-
40\left((f^{'})^3-(g^{'})^3\right)=0
\end{array}
\label{dfg}
\end{equation}
(here we used the identities $\partial_x=f^{'}\partial_f, ~
\partial_y=g^{'}\partial_g$). Applying to (\ref{dfg})
operator
$\partial^6/ \partial f^3\partial g^3$ we arrive at
$$
\frac{d^6((f^{'})^3)}{df^6}-\frac{d^6((g^{'})^3)}{dg^6}=0
$$
implying that $(f^{'})^3$ and $(g{'})^3$ are polynomials of the 6-th order in
$f$ and $g$, respectively. Coefficients of these polynomials are not
independent and can be fixed upon substitution in (\ref{dfg}):
\begin{equation}
(f^{'})^3=P(f), ~~~~
(g^{'})^3=P(-g)
\label{fg}
\end{equation}
where $P$ is an arbitrary polynomial of the 6th order.
Calculations  presented here follow \cite{Finikov37}, p.66-69.
Formulae (\ref{fg}) reflect the uniformizability of Kummer's quartics by
theta functions of genus 2. Since equations (\ref{3}) are invariant
under the transformation $(\beta, V, W)\to (-\beta, V, W)$,
the class of Kummer's quartics is self-dual.

Quartics of Kummer constitute  a subclass of

\bigskip

{\bf Isothermally asymptotic surfaces possessing a 3-parameter family
of projective applicabilities}
which are characterized by a condition
\begin{equation}
(\ln \beta^2)_{xy}=c\beta^2
\label{c}
\end{equation}
for some constant $c$. Quartics of Kummer correspond to $c=\frac{8}{9}$.
Condition (\ref{c}) means that the projective metric $2\beta^2~dxdy$ has
constant Gaussian curvature $K=-c$. To investigate equations (\ref{smVN})
with the additional constraint (\ref{c}) we introduce the anzatz
\begin{equation}
V=(\ln \beta)_{xx}+\frac{1}{2}(\ln \beta)_x^2+A, ~~~
W=(\ln \beta)_{yy}+\frac{1}{2}(\ln \beta)_y^2+B,
\label{AB}
\end{equation}
which implies upon substitution in (\ref{smVN}) the following equations
for $A, B$:
\begin{equation}
\begin{array}{c}
A_y=\frac{3}{2}(1-c)(\beta^2)_x, ~~~~ B_x=\frac{3}{2}(1-c)(\beta^2)_y, \\
\ \\
(\beta^2)_xA+\beta^2A_x=(\beta^2)_yB+\beta^2B_y.
\end{array}
\label{dAB}
\end{equation}
For $c=1$ equations (\ref{dAB}) are satisfied if $A=B=0$. The corresponding
surfaces are improper affine spheres; they will be discussed below. Here we
consider the case $c\ne 1$. Introducing $F$ by the formulae
\begin{equation}
A_x=-A(\ln \beta^2)_x+F, ~~~~ B_y=-B(\ln \beta^2)_y+F
\label{F}
\end{equation}
and writing down the compatibility conditions of (\ref{F}) with
$(\ref{dAB})_1$, $(\ref{dAB})_2$, we arrive at the equations for $F$
\begin{equation}
\begin{array}{c}
F_x=2cB\beta^2+\frac{3}{2}(1-c)\frac{1}{\beta^2}(\beta^2(\beta^2)_y)_y \\
\ \\
F_y=2cA\beta^2+\frac{3}{2}(1-c)\frac{1}{\beta^2}(\beta^2(\beta^2)_x)_x \\
\end{array}
\label{dF}
\end{equation}
the compatibility conditions of which coincide with (\ref{p}). Inserting in
(\ref{p}) the general solution of the Liouville equation (\ref{c})
\begin{equation}
\beta^2=\frac{1}{c}\frac{f^{'}g^{'}}{(f+g)^2}
\label{pc}
\end{equation}
we end up with the same $f, g$ as in (\ref{fg}). For any $\beta$ given
by (\ref{pc}), (\ref{fg}) the functions
$V$ and $W$ can be recowered from (\ref{AB}) where $A, B, F$
satisfy the compatible system (\ref{dAB}), (\ref{F}), (\ref{dF}).
Thus for any such $\beta$ there exists a 3-parameter family of surfaces which
have the same metric $2\beta^2~dxdy$ and the same cubic form
$\beta (dx^3+dy^3)$ and are not projectively equivalent. In general, two
projectively different surfaces in $P^3$ having the same projective metric
and the same cubic form  in a common asymptotic
parametrization $x, y$ are called projectively applicable. One can show
that only for isothermally asymptotic surfaces $M^2$ satisfying
(\ref{c}) does there exist a 3-parameter family of projectively different
surfaces which are all projectively applicable to $M^2$
(the value 3 is the maximal possible).
Isothermally asymptotic surfaces possessing only one-parameter families of
projective applicabilities have been investigated in \cite{Lingenberg2}.

The whole class of projectively applicable surfaces (not necessarily
isothermally asymptotic) is discussed sect.4.

\bigskip

{\bf Affine spheres} constitute an important subclass of isothermally
asymptotic surfaces specified by the following reduction in (\ref{smVN}):
\begin{equation}
V=\frac{\beta_{xx}}{\beta}-\frac{1}{2}\left(\frac{\beta_x}{\beta}\right)^2,
~~~~~
W=\frac{\beta_{yy}}{\beta}-\frac{1}{2}\left(\frac{\beta_y}{\beta}\right)^2.
\label{VW}
\end{equation}
After this ansatz the first equation in (\ref{smVN}) will be satisfied
identically while the last two imply the Tzitzeica equation for $\beta$:
\begin{equation}
(\ln \beta)_{xy}=\beta^2+\frac{c}{\beta}, ~~~ c=const.
\label{Tzitzeica}
\end{equation}
The cases $c\ne 0$ and $c=0$ correspond to proper and improper affine spheres,
respectively. In view of (\ref{klab}) equations (\ref{VW}) and
(\ref{Tzitzeica}) are equivalent to
$$
a=b=0, ~~~~  k=l=-\frac{c}{\beta}.
$$
Using equations (\ref{Wilczynski}) it is easy to check that the point
$$
{\sqrt {\beta}}~{\mbox{\boldmath $\eta$}}+\frac{c}{2{\sqrt {\beta}}}~{\bf r}
$$
is independent of $x, y$. Geometrically, this means that all normals of
Wilczynski intersect in one fixed point.
This property can be viewed as a projective-invariant
definition of affine spheres (we emphasize that the normals of Wilczynski do not
coincide in general with Blaschke's normals in affine differential geometry).

With $V, W$ given by (\ref{VW}) equations (\ref{r}) possess a
particular solution $r^0=\sqrt \beta$. Upon the substitution
$R=r/r^0$ equations (\ref{affine}) assume the form
$$
\begin{array}{c}
R_{xx}=\beta R_y-\frac{\beta_x}{\beta} R_x \\
\ \\
R_{yy}=\beta R_x-\frac{\beta_y}{\beta} R_y \\
\end{array}
$$
which become the familiar equations for the radius-vector of affine spheres
after adding the compatible equation
$$
R_{xy}=-\frac{c}{\beta}R.
$$
The fact that Tzitzeica's equation solves the stationary mVN equation is also
reflected in the following nonlocal representation of the mVN
equation:
$$
\beta_t=\left(\frac{1}{\beta}\partial_x \beta^2 \partial_y^{-1}\frac{1}{\beta}
\partial_x - \frac{1}{\beta}\partial_y \beta^2
\partial_x^{-1}\frac{1}{\beta}\partial_y\right)
  \left(\beta(\ln \beta)_{xy}-\beta^3\right).
$$
Indeed, the condition $\beta (\ln \beta)_{xy}-\beta^3=c=const$ is equivalent to
\ref{Tzitzeica}).

\bigskip

We refer to \cite{Lane}, \cite{Bol1},
\cite{Finikov37}, \cite{Finikov50}, \cite{Fer9} for the further discussion
of isothermally asymptotic surfaces.

In the recent paper \cite{KonPin}
Konopelchenko and Pinkall introduced integrable dynamics of surfaces in
3-space governed by the Veselov-Novikov equation. One can show that
isothermally asymptotic surfaces can be interpreted as the stationary points
of this evolution. So it is not surprising that they are described
by the stationary mVN equation (which, as it has been demonstrated in
\cite{Fer9}, is equivalent to the stationary VN).

\section{Projectively applicable surfaces}

Generically, the quadratic form
$$
\beta \gamma\,dxdy
$$
and the cubic form
$$
\beta \,dx^3+\gamma \,dy^3
$$
define a surface $M^2\subset P^3$ uniquely up to projective equivalence.
There exists, however, a class of the so-called projectively applicable
surfaces for which this is not the case.
Projective deformations of surfaces have been extensively
investigated by Cartan \cite{Cartan1}.
Analytically, projective applicability
means that for a given $\beta, \gamma$ equations (\ref{GC1}) do not specify
$V$ and $W$ uniquely (equivalently, equations (\ref{GC2})$_3$ --
(\ref{GC2})$_5$
do not uniquely specify $a$ and $b$). As shown in \cite[p.\ 62]{Finikov37},
projectively applicable surfaces naturally fall into two different classes:
the so-called "surfaces $R_0$" and "surfaces $R$".

\bigskip

{\bf Surfaces ${\bf R_0}$} are specified by the condition $(\ln \beta)_{xy}=0$
(or $(\ln \gamma)_{xy}=0$),
which can be reduced to $\beta=1$ (or $\gamma=1$) after a suitable choice of
$x, y$. The substitution of $\beta=1$ in (\ref{GC1}) results in the system
\begin{equation}
\begin{array}{c}
\gamma_{xxx}-2\gamma_xV-\gamma V_x+W_y=0\\
\ \\
W_x=\gamma_y \\
\ \\
V_y=2 \gamma_x.
\end{array}
\label{R_0}
\end{equation}
Since $W$ enters  (\ref{R_0}) only through it's derivatives, it is defined up
to an
additive constant $W\to W+\lambda$, providing thus a 1-parameter family
of projectively nonequivalent surfaces having the same $\beta$ and $\gamma$.
The radius-vectors of these surfaces satisfy the linear system
$$
\begin{array}{c}
{\bf r}_{xx}={\bf r}_y+\frac{1}{2} V {\bf r} \\
\ \\
{\bf r}_{yy}=\gamma {\bf r}_x+\frac{1}{2}(W-\gamma_x+\lambda) {\bf r}
\end{array}
$$
containing a "spectral parameter" $\lambda$. Rewriting this linear
system in the form
\begin{equation}
\begin{array}{c}
{\bf r}_{y}={\bf r}_{xx}-\frac{1}{2} V {\bf r} \\
\ \\
{\bf r}_{xxxx}=V {\bf r}_{xx}+(V_x+\gamma){\bf r}_x+
\frac{1}{2}(V_{xx}+\gamma_x-\frac{1}{2}V^2+W+\lambda) {\bf r}
\end{array}
\label{statKP}
\end{equation}
we immediately recognise in (\ref{statKP})$_1$ the nonstationary Schr\"odinger
equation which plays a role of the Lax operator of the (2+1)-dimensional
integrable Kadomtsev-Petviashvili (KP) hierarchy. Meanwhile, (\ref{statKP})$_2$
coincides with the stationary limit of the linear problem generating the
second flow in the KP hierarchy.

We recall that the second flow in the KP hierarchy is associated with the
linear system
\begin{equation}
\begin{array}{c}
{\bf r}_{y}={\bf r}_{xx}-\frac{1}{2} V {\bf r} \\
\ \\
{\bf r}_t={\bf r}_{xxxx}-V {\bf r}_{xx}-(V_x+\gamma){\bf r}_x-
\frac{1}{2}(V_{xx}+\gamma_x-\frac{1}{2}V^2+W) {\bf r}
\end{array}
\label{KP}
\end{equation}
which is related to (\ref{statKP}) via a formal transformation
$\frac{1}{2}\lambda {\bf r} \to \partial_t {\bf r}$. The compatibility
conditions of (\ref{KP}) result in the "second" KP flow:
$$
\begin{array}{c}
V_t=\gamma_{xxx}-2\gamma_xV-\gamma V_x+W_y\\
\ \\
W_x=\gamma_y \\
\ \\
V_y=2 \gamma_x.
\end{array}
$$
Thus, equations (\ref{R_0}) governing $R_0$-surfaces coincide with the
stationary limit of the second flow in the KP hierarchy. In a sense, they
can be called the higher Boussinesq equations.

Geometric properties of surfaces $R_0$ have been investigated, e.g., in
\cite{Bol1}, \cite{Mihailescu}.

\bigskip

{\bf Surfaces R} are specified by the condition $\beta_y=\gamma_x$. In this
case
$V$ and $W$ are defined by (\ref{GC1}) up to additive constant
$W\to W+\lambda, ~ V\to V+\lambda$, so that we again arrive at a 1-parameter
family
of projectively nonequivalent surfaces having the same $\beta, \gamma$.
Their radius-vectors satisfy the linear system
$$
\begin{array}{c}
{\bf r}_{xx}=\beta {\bf r}_y+\frac{1}{2}(V-\beta_y+\lambda) {\bf r} \\
\ \\
{\bf r}_{yy}=\gamma {\bf r}_x+\frac{1}{2}(W-\gamma_x+\lambda) {\bf r}
\end{array}
$$
containing a "spectral parameter" $\lambda$. Rewriting it in the form
$$
\begin{array}{c}
{\bf r}_{xx}-{\bf r}_{yy}=\beta {\bf r}_y-\gamma {\bf r}_x+
\frac{1}{2}(V-W) {\bf r} \\
\ \\
{\bf r}_{xx}+{\bf r}_{yy}-\beta {\bf r}_y-\gamma {\bf r}_x-
\frac{1}{2}(V+W-\beta_y-\gamma_x) {\bf r}=\lambda {\bf r},
\end{array}
$$
we recognize the stationary DS-type linear problem
(the relationship of surfaces R to the stationary DS hierarchy has been pointed
out in \cite{Kon2}). Replacing $\lambda {\bf r}$ by $\partial_t {\bf r}$
we arrive at the (2+1)-dimensional DS linear problem
$$
\begin{array}{c}
{\bf r}_{xx}-{\bf r}_{yy}=\beta {\bf r}_y-\gamma {\bf r}_x+
\frac{1}{2}(V-W) {\bf r} \\
\ \\
{\bf r}_t={\bf r}_{xx}+{\bf r}_{yy}-\beta {\bf r}_y-\gamma {\bf r}_x-
\frac{1}{2}(V+W-\beta_y-\gamma_x) {\bf r},
\end{array}
$$
with the compatibility conditions
\begin{equation}
\begin{array}{c}
V_t+\beta_{yyy}-2\beta_yW-\beta W_y=
W_t+\gamma_{xxx}-2\gamma_xV-\gamma V_x\\
\ \\
\gamma_t=W_x-2\gamma \beta_y-\beta \gamma_y \\
\ \\
\beta_t=V_y-2\beta \gamma_x-\gamma \beta_x \\
\ \\
\beta_y=\gamma_x.
\end{array}
\label{R}
\end{equation}
In the stationary limit they reduce to the equations governing surfaces R.
Introducing
$$
\beta=\varphi_x, ~~~ \gamma=\varphi_y, ~~~
V=\frac{1}{2}\varphi_x^2 +\theta_{xx}, ~~~
W=\frac{1}{2}\varphi_y^2 +\theta_{yy},
$$
we can rewrite (\ref{R}) in the form
$$
\begin{array}{c}
\varphi_t=\theta_{xy}-\varphi_x\varphi_y, \\
\ \\
\theta_t=\varphi_{xy}-P, \\
\ \\
\Box P=(\varphi_x \Box \theta)_y+(\varphi_y \Box \theta)_x,
\end{array}
$$
where $\Box =\partial_x^2-\partial_y^2$.

The coordinate net $\xi=x+y, ~ \eta=x-y$ on a surface R is conjugate and has
an
important geometric property, namely, that both congruences of lines, tangent
to the $\xi$- and $\eta$-coordinate curves, form a W-congruences.
Conjugate nets with this property are called the R-nets; they have been
introduced by
Demoulin and Tzitzeica in a series of papers \cite{Demoulin1}, \cite{Tzi11}.
Thus, surfaces R can be equivalently characterized by the existence of an
R-net.
B\"acklund transformation for surfaces R was established by Jonas \cite{Jonas1}
and subsequently discussed in \cite{Finikov37}, \cite{Finikov50},
\cite{Eisenhart}. Some further properties of surfaces $R$ have been studied
in \cite{Gambier}, \cite{Godeaux1}. Multidimensional analogs of surfaces $R$
were investigated in a series of papers \cite{Korovin}, \cite{Geidelman},
\cite{Finikov57}, \cite{Shulman2}.

The most important examples of  surfaces R include:

\bigskip

{\bf Projective transforms of surfaces with constant Gaussian curvature
$K=-1$ in the Euclidean 3-space} which correspond to
$$
\begin{array}{c}
\beta=-\frac{\varphi_x}{\sin \varphi}, ~~~
\gamma=-\frac{\varphi_y}{\sin \varphi}, \\
\ \\
V=1+\frac{1}{2}\beta^2\cos^2\varphi+(\beta\cos\varphi)_x, ~~~~
W=1+\frac{1}{2}\gamma^2\cos^2\varphi+(\gamma\cos\varphi)_y
\end{array}
$$
where $\varphi$ satisfies the Sine-Gordon equation
$\varphi _{xy}=-\sin\varphi$. One can check directly that with this anzatz
equations (\ref{GC1}) as well as $\beta_y=\gamma_x$ are satisfied identically.
Moreover, equations (\ref{r}) possess a particular solution
$r^0=\frac{1}{\sqrt {\sin \varphi}}$, so that in the affine gauge
${\bf R}={\bf r}/r^0={\sqrt {\sin \varphi }}\ {\bf r}$
equations (\ref{r}) transform to
$$
\begin{array}{c}
{\bf R}_{xx}=-\frac{\varphi_x}{\sin \varphi}{\bf R}_y+
\frac{\varphi_x \cos \varphi}{\sin \varphi}{\bf R}_x \\
\ \\
{\bf R}_{yy}=-\frac{\varphi_y}{\sin \varphi}{\bf R}_x+
\frac{\varphi_y \cos \varphi}{\sin \varphi}{\bf R}_y.
\end{array}
$$
Supplementing these equations with
$$
\begin{array}{c}
{\bf R}_{xy}=\sin \varphi \ {\bf n} \\
\ \\
{\bf n}_x=\frac{1}{\sin \varphi}({\bf R}_y-\cos \varphi {\bf R}_x) \\
\ \\
{\bf n}_y=\frac{1}{\sin \varphi}({\bf R}_x-\cos \varphi {\bf R}_y)
\end{array}
$$
(which are all mutually compatible in view of $\varphi_{xy}=-\sin \varphi$,
and, moreover, possess a specialization $({\bf n}, {\bf n})=
({\bf R}_x, {\bf R}_x)=({\bf R}_y, {\bf R}_y)=1, ~ ({\bf R}_x, {\bf R}_y)=
\cos \varphi)$, we immediately recognize the equations governing the
radius-vector ${\bf R}$ and the unit normal ${\bf n}$
of a surface with constant Gaussian
curvature $K=-1$ parametrized by asymptotic coordinates $x, y$. It's first and
second
fundamental forms are given by the formulae
$$
\begin{array}{c}
dx^2+2\cos \varphi  \ dxdy+dy^2\\
\ \\
2\sin\varphi \ dxdy,
\end{array}
$$
respectively. In this case the conjugate net R is the net of curvature lines.

Closely related examples of surfaces R are focal surfaces of the congruence
of normals of a surface with $K=-1$.

\bigskip

{\bf Projective transforms of surfaces with constant Gaussian curvature $K=1$
in the Lorentzian 3-space} are characterized by exactly the same formulae as
surfaces with $K=-1$ in the Euclidean space with the interchange
$\sin \to \sinh, \ \cos \to \cosh$ in all places where they occur.
These surfaces are described by the
equation $\varphi_{xy}=-\sinh \varphi$ and have the fundamental forms
$$
\begin{array}{c}
dx^2+2\cosh \varphi  \ dxdy+dy^2\\
\ \\
2\sinh \varphi \ dxdy,
\end{array}
$$
respectively. Note that the first fundamental form is now indefinite.

\bigskip

{\bf Surfaces with an R-net of period 4} are characterized by
$$
\begin{array}{c}
\beta=\frac{1}{2}\varphi_x, ~~~ \gamma=\frac{1}{2}\varphi_y \\
\ \\
V=\frac{1}{4}\varphi_x^2+\frac{1}{8}\varphi_y^2+\sinh \varphi, ~~~~
W=\frac{1}{4}\varphi_y^2+\frac{1}{8}\varphi_x^2-\sinh \varphi \\
\end{array}
$$
where $\varphi_{xx}-\varphi_{yy}=4 \cosh \varphi$. With these $\beta, \gamma,
V, W$
equations (\ref{GC1}) and $\beta_y=\gamma_x$ are satisfied identically.
One can show that under  Laplace transformations the R-net
$\xi=x+y, \ \eta=x-y$ on these surfaces is periodic with period four.
Moreover, all surfaces carrying an R-net of period four
can be obtained in this way.

\section{Surfaces of Jonas}

are characterized by a condition
$$
\beta_x=\gamma_y
$$
which is equivalent to the requirement that the conjugate net
$\xi=x+y, \ \eta=x-y$ on the surface $M^2$ has equal point and tangential
Laplace invariants (tangential Laplace invariants are the Laplace invariants
of the conjugate net $\xi, \eta$ on the dual surface).
This class of surfaces has been introduced by Jonas in \cite{Jonas2} as a
counterpart of surfaces R and subsequently discussed in
\cite{Cartan}, \cite{Gambier}, \cite{Pa}  as well as in the textbooks
cited above.
Conjugate nets with equal point and tangential Laplace invariants are
called the nets of Jonas.
The condition $\beta_x=\gamma_y$ allows an introduction of a spectral parameter
$\lambda$ in the linear system  (\ref{UAPVBQ}). For that purpose we first
define
$H$ and $K$ by the formulae
\begin{equation}
\begin{array}{c}
H_x=-\gamma \ K +\lambda \ {\cal V}, ~~~~ H_y=-\beta \  K \\
\ \\
K_y=-\beta \ H +\lambda \ {\cal U}, ~~~~ K_x=-\gamma \  H
\end{array}
\label{HK}
\end{equation}
which are compatible in view of (\ref{UAPVBQ}) and $\beta_x=\gamma_y$.
The systems (\ref{UAPVBQ}) and (\ref{HK}) can be coupled in an $8\times 8$
linear system
\begin{equation}
\begin{array}{c}
\left(\begin{array}{c}
{\cal U}\\
{\cal A}\\
{\cal P}\\
{\cal V}\\
{\cal B}\\
{\cal Q}\\
H\\
K
\end{array}\right)_x=
\left(\begin{array}{cccccccc}
0 & 0 & 0 & \beta & 0 & 0 & 0 & 0\\
k & 0 & 0 & 0 & 0 & 0 & 0 & 0\\
0 & k & 0 & -\beta a & 0 & 0 & 0 & 0\\
0 & 0 & 0 & \frac{\gamma_x}{\gamma} & 1 & 0 & 0 & 0\\
0 & 0 & 0 & b & 0 & 1 & 0 & 0\\
-\beta a & 0 & \beta & 0 & b &-\frac{\gamma_x}{\gamma} & 1 & 0\\
0 & 0 & 0 & \lambda & 0 & 0 & 0 & -\gamma\\
0 & 0 & 0 & 0 & 0 & 0 & -\gamma & 0
\end{array}\right)
\left(\begin{array}{c}
{\cal U}\\
{\cal A}\\
{\cal P}\\
{\cal V}\\
{\cal B}\\
{\cal Q}\\
H\\
K
\end{array}\right)\\
\ \\
\left(\begin{array}{c}
{\cal U}\\
{\cal A}\\
{\cal P}\\
{\cal V}\\
{\cal B}\\
{\cal Q}\\
H\\
K
\end{array}\right)_y=
\left(\begin{array}{cccccccc}
\frac{\beta_y}{\beta} & 1 & 0 & 0 & 0 & 0 & 0 & 0\\
a & 0 & 1 & 0 & 0 & 0 & 0 & 0\\
0 & a & -\frac{\beta_y}{\beta} & -\gamma b & 0 & \gamma & 0 & 1\\
\gamma & 0 & 0 & 0 & 0 & 0 & 0 & 0\\
0 & 0 & 0 & l & 0 & 0 & 0 & 0\\
-\gamma b & 0 & 0 & 0 & l & 0 & 0 & 0\\
 0 & 0 & 0 & 0 & 0 & 0 & 0 & -\beta\\
\lambda & 0 & 0 & 0 & 0 & 0 & -\beta & 0
\end{array}\right)
\left(\begin{array}{c}
{\cal U}\\
{\cal A}\\
{\cal P}\\
{\cal V}\\
{\cal B}\\
{\cal Q}\\
H\\
K
\end{array}\right)
\end{array}
\label{8x8}
\end{equation}

The coupling is achieved by adding $K$ to ${\cal P}_y$ and
$H$ to ${\cal Q}_x$. Linear system (\ref{8x8}) reduces to (\ref{UAPVBQ})
under the reduction $\lambda=H=K=0$. This linear system is gauge equivalent
to the one used in \cite{Jonas2} for the construction of the B\"acklund
transformation for the Jonas surfaces.
With $\beta=\varphi_y, \ \gamma=\varphi_x$ equations (\ref{GC1})
assume the form
\begin{equation}
\begin{array}{c}
\varphi_{yyyy}-2 \ \varphi_{yy} \ W-\varphi_y \ W_y=
\varphi_{xxxx}-2 \ \varphi_{xx} \ V-\varphi_x \ V_x\\
\ \\
(W+\frac{1}{2}\varphi_y^2)_x=2 \ (\varphi_x \varphi_y)_y \\
\ \\
(V+\frac{1}{2}\varphi_x^2)_y=2 \ (\varphi_x \varphi_y)_x.
\end{array}
\label{Jonas}
\end{equation}
Rewriting equations (\ref{8x8}) in terms of ${\cal U},  {\cal V}$
(that is, expressing
${\cal A}, {\cal B}, {\cal P}, {\cal Q}, H, K$ through
${\cal U}, {\cal U} _y, {\cal U} _{yy}, {\cal U} _{yyy}$ and
${\cal V}, {\cal V} _x,
{\cal V} _{xx}, {\cal V} _{xxx}$), we arrive at
the equivalent spectral problem of the system (\ref{Jonas}):
\begin{equation}
\begin{array}{c}
{\cal U}_x=\varphi_y \ {\cal V} \\
{\cal V}_y=\varphi_x \ {\cal U}
\end{array}
\label{linJonas1}
\end{equation}

\begin{equation}
\begin{array}{c}
\lambda \ {\cal U}={\cal U}_{xxxx}+{\cal U}_{yyyy}-(\varphi_y^2+2W)\ {\cal
U}_{yy}-
4 \ \varphi_{xy} \ {\cal V}_{xx} \\
\ \\
+(\varphi_y \ \varphi_{yy} - 3\ W_y)\ {\cal U}_y-
(2\ \varphi_{xxy}+\varphi_y \ \varphi_x^2+2 \ V \varphi_y)\ {\cal V}_x+
m \ {\cal U} +n \ {\cal V}, \\
\ \\
\ \\
\lambda \ {\cal V}={\cal V}_{xxxx}+{\cal V}_{yyyy}-(\varphi_x^2+2V)\ {\cal
V}_{xx}-
4 \ \varphi_{xy} \ {\cal U}_{yy} \\
\ \\
+(\varphi_x \ \varphi_{xx} - 3\ V_x)\ {\cal V}_x-
(2\ \varphi_{xyy}+\varphi_x \ \varphi_y^2+2 \ W \varphi_x)\ {\cal U}_y+
\tilde m \ {\cal V} +\tilde n \ {\cal U},
\end{array}
\label{linJonas2}
\end{equation}
where $\lambda = const$ and $m, \ \tilde m, \ n, \ \tilde n$ are given by the
formulae
$$
\begin{array}{c}
m=-W_{yy}-\varphi_y \ \varphi_{yyy} - 2 \ \varphi_x \ \varphi_{xxx}
+\varphi_{xx}^2+2 \ V \ \varphi_x^2+2 \ W \ \varphi_y^2 \\
\ \\
\tilde m=-V_{xx}-\varphi_x \ \varphi_{xxx} - 2 \ \varphi_y \ \varphi_{yyy}
+\varphi_{yy}^2+2 \ V \ \varphi_x^2+2 \ W \ \varphi_y^2 \\
\ \\
n=-2 \ \varphi_{xxxy}+\varphi_x^2 \ \varphi_{xy}+3 \ \varphi_x \ \varphi_y \
\varphi_{xx}+ 2 \ V \ \varphi_{xy} - \varphi_y \ V_x \\
\ \\
\tilde n=-2 \ \varphi_{xyyy}+\varphi_y^2 \ \varphi_{xy}+3 \ \varphi_x \
\varphi_y \
\varphi_{yy}+ 2 \ W \ \varphi_{xy} - \varphi_x \ W_y.
\end{array}
$$
The structure of linear system (\ref{linJonas1}),  (\ref{linJonas2})
clearly suggests that equations  (\ref{Jonas}) are related to
the stationary limit of the fourth order
flow in the DS hierarchy. Indeed, the fourth order
flow in the DS hierarchy is generated by the 2-dimensional Dirac operator
\begin{equation}
\begin{array}{c}
{\cal U}_x=\beta \ {\cal V}
 \\
{\cal V}_y=\gamma \ {\cal U}
\end{array}
\label{DS4}
\end{equation}
where the time evolution of ${\cal U}$ and ${\cal V}$ is specified by
\begin{equation}
\begin{array}{c}
{\cal U}_t={\cal U}_{xxxx}+{\cal U}_{yyyy}+...\\
\ \\
{\cal V}_t={\cal V}_{xxxx}+{\cal V}_{yyyy}+....
\end{array}
\label{DS4t}
\end{equation}
The linear problem corresponding to the stationary flow can be obtained
by a formal substitution ${\cal U}_t\to \lambda {\cal U}, \
{\cal V}_t\to \lambda {\cal V}$ where
$\lambda$ plays a role of spectral parameter. This linear
problem reduces to (\ref{linJonas1}), (\ref{linJonas2})
 after the additional reduction
$\beta_x=\gamma_y$. We emphasize that the reduction $\beta_x=\gamma_y$
is compatible not with the  (2+1)-dimensional flow generated by
(\ref{DS4}), (\ref{DS4t}), rather than only with its stationary limit. This
reduction
is compatible only with the odd order flows in the DS hierarchy.

Particular examples of the Jonas surfaces are provided by

\bigskip
{\bf Projective transforms of minimal surfaces in the Euclidean 3-space}
corresponding to
$$
\begin{array}{c}
\beta=\varphi_y, ~~~
\gamma=\varphi_x, \\
\ \\
V=\frac{1}{2}\varphi_x^2 - \varphi_y^2 - \exp {2\varphi}, ~~~~
W=\frac{1}{2}\varphi_y^2 - \varphi_x^2 - \exp {2\varphi}
\end{array}
$$
where $\varphi$ satisfies the Liouville equation
$\varphi _{xx}+\varphi _{yy}=-\exp {2\varphi}$.
One can check directly that with this anzatz
equations (\ref{GC1}) as well as $\beta_x=\gamma_y$ are satisfied identically.
Moreover, equations (\ref{r}) possess a particular solution
$r^0=\exp \frac{\varphi}{2}$, so that in the affine gauge
${\bf R}={\bf r}/r^0= \exp (-\frac{\varphi}{2}){\bf r}$
equations (\ref{r}) transform to
$$
\begin{array}{c}
{\bf R}_{xx}=\varphi_y \ {\bf R}_y-\varphi_x \ {\bf R}_x \\
{\bf R}_{yy}=\varphi_x \ {\bf R}_x-\varphi_y \ {\bf R}_y.
\end{array}
$$
Supplementing these equations with
$$
\begin{array}{c}
{\bf R}_{xy}=-\varphi_y \ {\bf R}_x-\varphi_x \ {\bf R}_y - {\bf n} \\
\ \\
{\bf n}_x=\exp {2\varphi} \ {\bf R}_y \\
\ \\
{\bf n}_y=\exp {2 \varphi} \ {\bf R}_x
\end{array}
$$
which are all mutually compatible in view of
$\varphi _{xx}+\varphi _{yy}=-\exp {2\varphi}$,
and, moreover, possess a specialization $({\bf n}, {\bf n})=1, \
({\bf R}_x, {\bf R}_x)=({\bf R}_y, {\bf R}_y)=\exp {(-2 \varphi)}, ~
({\bf R}_x, {\bf R}_y)=0$, we immediately recognize the equations governing the
radius-vector ${\bf R}$ and the unit normal ${\bf n}$
of a minimal surface parametrized by asymptotic coordinates $x, y$.
It's first and second fundamental forms are given by the formulae
$$
\begin{array}{c}
\exp {(-2\varphi)}(dx^2+dy^2)\\
\ \\
-2 \ dxdy,
\end{array}
$$
respectively. For minimal surfaces
the conjugate net of Jonas is the net of curvature lines.

In \cite{Jonas2} there was established a one-to-one correspondence between
surfaces of Jonas and pairs of surfaces in the Euclidean 3-space such that:
\par
\noindent -- both surfaces are in isometric correspondence;
\par
\noindent -- this correspondence preserves a congugate net with equal point
Laplace invariants.
\par
\noindent In a sense, this construction is an analog of the correspondence
between isothermic surfaces and Bonnet pairs (Bonnet pairs are pairs of
surfaces in isometric correspondence such that the mean curvature
in the corresponding points is the same).

\section{Projectively minimal surfaces}

These surfaces
are the extremals of  the projective area functional
\begin{equation}
\int \int \beta\gamma\,dxdy.
\label{int}
\end{equation}
The Euler-Lagrange equations for the functional
(\ref{int}) have been derived in \cite{Thomsen1}, \cite{Sasaki} and
in our notation adopt the form
$$
\begin{array}{c}
\beta_{yyy}-2\beta_yW-\beta W_y=0, ~~~~
\gamma_{xxx}-2\gamma_xV-\gamma V_x=0\\
\ \\
W_x=2\gamma \beta_y+\beta \gamma_y \\
\ \\
V_y=2\beta \gamma_x+\gamma \beta_x
\end{array}
$$
which results after equating to zero both sides of (\ref{GC1})$_1$.
Exploiting the obvious symmetry $\beta \to \lambda \beta, \
\gamma \to \frac{1}{\lambda} \gamma$, we can introduce the spectral parameter
in the
linear system (\ref{r}):
$$
\begin{array}{c}
{\bf r}_{xx}=\lambda \ \beta \ {\bf r}_y+
\frac{1}{2}(V-\lambda \ \beta_y) \ {\bf r} \\
\ \\
{\bf r}_{yy}=\frac{1}{\lambda} \ \gamma \ {\bf r}_x+
\frac{1}{2}(W-\frac{1}{\lambda} \ \gamma_x) \ {\bf r}.
\end{array}
$$
In terms of (\ref{GC2}) the Euler-Lagrange equations reduce to
$$
\beta a_y+2a\beta_y=0, ~~~ \gamma b_x+2b\gamma_x=0,
$$
which result after equating to zero both sides of equation
(\ref{GC2})$_5$. In the Pl\"ucker coordinates the
spectral problem assumes the form
\begin{equation}
\begin{array}{c}
\left(\begin{array}{c}
{\cal U}\\
{\cal A}\\
{\cal P}\\
{\cal V}\\
{\cal B}\\
{\cal Q}
\end{array}\right)_x=
\left(\begin{array}{cccccc}
0 & 0 & 0 & \lambda \beta & 0 & 0\\
k & 0 & 0 & 0 & 0 & 0\\
0 & k & 0 & -\lambda \beta a & 0 & 0\\
0 & 0 & 0 & \frac{\gamma_x}{\gamma} & 1 & 0\\
0 & 0 & 0 & b & 0 & 1\\
-\lambda \beta a & 0 & \lambda \beta & 0 & b &-\frac{\gamma_x}{\gamma}
\end{array}\right)
\left(\begin{array}{c}
{\cal U}\\
{\cal A}\\
{\cal P}\\
{\cal V}\\
{\cal B}\\
{\cal Q}
\end{array}\right)\\
\ \\
\left(\begin{array}{c}
{\cal U}\\
{\cal A}\\
{\cal P}\\
{\cal V}\\
{\cal B}\\
{\cal Q}
\end{array}\right)_y=
\left(\begin{array}{cccccc}
\frac{\beta_y}{\beta} & 1 & 0 & 0 & 0 & 0\\
a & 0 & 1 & 0 & 0 & 0\\
0 & a & -\frac{\beta_y}{\beta} & -\frac{\gamma b}{\lambda} & 0 & \frac
{\gamma}{\lambda}\\
\frac{\gamma}{\lambda} & 0 & 0 & 0 & 0 & 0\\
0 & 0 & 0 & l & 0 & 0\\
-\frac{\gamma b}{\lambda} & 0 & 0 & 0 & l & 0
\end{array}\right)
\left(\begin{array}{c}
{\cal U}\\
{\cal A}\\
{\cal P}\\
{\cal V}\\
{\cal B}\\
{\cal Q}
\end{array}\right)
\end{array}
\label{Laxprojmin}
\end{equation}
Setting
$$
a=\frac{\varphi(x)}{\beta^2}, ~~~ b=\frac{\psi(y)}{\gamma^2},
$$
 we have three cases to distinguish:

{\bf General case}. Both $\varphi(x)$ and $\psi(y)$ are nonzero.
In this
case, we can always normalize $\varphi(x), \psi(y)$ to $\pm 1$ by means of
transformations (\ref{new}). Let us assume, for instance, that
$\varphi(x)= \psi(y)=1$. With this normalization equations (\ref{GC2}) assume
the form
\begin{equation}
\begin{array}{c}
(\ln \beta)_{xy}=\beta \gamma -k, ~~~~ (\ln \gamma)_{xy}=\beta \gamma -l,\\
\ \\
(\beta k)_y+2\frac{\beta_x}{\beta^2}=0, ~~~~
(\gamma l)_x+2\frac{\gamma_y}{\gamma^2}=0.
\end{array}
\label{projmin}
\end{equation}

{\bf Surfaces of Godeaux-Rozet} \cite[p.\ 318]{Bol1}. In this case,
$\varphi =0$, and hence $a=0$,
while $\psi$ is nonzero and may be normalized to $\pm 1$.
Here, we assume that $\psi=1$. Inserting this ansatz in (\ref{GC2}) we obtain
$$
k=\frac{s(x)}{\beta}
$$
Hence, if $s(x)$ is nonzero, it may be reduced to $-1$
by means of (\ref{new}) so that the resulting equations take the form
\begin{equation}
\begin{array}{c}
(\ln \beta)_{xy}=\beta \gamma +\frac{1}{\beta},
~~~~ (\ln \gamma)_{xy}=\beta \gamma -l,\\
\ \\
(\gamma l)_x+2\frac{\gamma_y}{\gamma^2}=0.
\end{array}
\label{Godeaux-Rozet}
\end{equation}

{\bf Surfaces of Demoulin}. In this case,
both $\varphi$ and $\psi$ are zero and hence $a=b=0$, so that
$$
k=\frac{s(x)}{\beta}, ~~~ l=\frac{t(y)}{\gamma}.
$$
Once again, the analysis falls into three subcases depending on whether
$s,t$ are zero or not. In the generic situation
$s\ne 0, \  t\ne 0$ both $s$ and $t$ may be normalized to $-1$  and the
resulting equations reduce to the coupled Tzitzeica system
\begin{equation}
\begin{array}{c}
(\ln \beta)_{xy}=\beta \gamma +\frac{1}{\beta}, \\
\ \\
(\ln \gamma)_{xy}=\beta \gamma +\frac{1}{\gamma},
\end{array}
\label{Demoulin}
\end{equation}
In this form, the equations governing Demoulin surfaces have been set down in
\cite[p.\ 51]{Finikov37}. The same system has been presented in
\cite{Mikhailov}
as a reduction of the two-dimensional Toda lattice. Inserting
$a=b=0, \ k=-\frac{1}{\beta}, \ l=-\frac{1}{\gamma}$ in  equations
(\ref{Laxprojmin}) and rewriting them in terms of $\cal A$ and $\cal B$, we
arrive at a more convenient second-order linear problem \cite{Fer8}
\begin{equation}
\begin{array}{c}
{\cal B}_{xy}=-\frac{1}{\gamma}\ {\cal B} ~~~~~~~~~~~~~~~~~~~
{\cal A}_{xy}=-\frac{1}{\beta}\ {\cal A} \\
\ \\
{\cal B}_{xx}=\lambda \ \beta \ {\cal A}_y-\frac{\gamma_x}{\gamma}\ {\cal B}_x
~~~~~~~~
{\cal A}_{xx}=\lambda \ \gamma \ {\cal B}_y-\frac{\beta_x}{\beta}\ {\cal A}_x
\\
\ \\
{\cal B}_{yy}=\frac{\beta}{\lambda} \ {\cal A}_x-\frac{\gamma_y}{\gamma}\
{\cal B}_y ~~~~~~~~~~~~
{\cal A}_{yy}=\frac{\gamma}{\lambda} \ {\cal B}_x-\frac{\beta_y}{\beta}\
{\cal A}_y.
\end{array}
\label{second-order}
\end{equation}

In \cite{Demoulin} Demoulin
established in a purely geometric manner the existence of B\"acklund
transformations for Godeaux-Rozet and Demoulin surfaces and proved the
corresponding
permutability theorems. Apparently, Demoulin did not formulate his results
in terms of analytic expressions. In \cite{Fer8}, a Toda lattice connection
and the second-order linear problem (\ref{second-order}) is used to derive
explicitly a B\"acklund transformation for
Demoulin surfaces.

{\bf Remark.} The specialization $\beta=\gamma$
reduces (\ref{Demoulin}) to the Tzitzeica
equation
$$
(\ln \beta)_{xy}=\beta^2 +\frac{1}{\beta},
$$
which governs affine spheres in affine differential
geometry \cite{Tzi10}. Moreover, with ${\cal A}={\cal B}$ the linear system
(\ref{second-order}) becomes the standard Lax pair for the Tzitzeica equation.
Geometrically this means that affine spheres lie in the
intersection of two different integrable classes of projective
surfaces, namely isothermally asymptotic and projectively minimal surfaces.

\medskip
Projectively minimal, Godeaux-Rozet and Demoulin surfaces
also arise in the theory of envelopes of Lie quadrics associated with the
surface $M^2$. For brevity, we only recall the necessary definitions. The
details can be found in \cite{Bol1}, \cite{Finikov37}, \cite{Lane}, etc.
Thus, let us consider a point $p^0$ on the
surface $M^2$ and the $x$-asymptotic line passing through $p^0$. Let us take
three additional points $p^i,\,i=1,2,3$ on this asymptotic line close to $p^0$
and draw three $y$-asymptotic lines $\gamma^i$ passing through $p^i$.
The three
straight lines which are tangential to $\gamma^i$ and pass through the points
$p^i$ uniquely define a quadric ${\bf Q}$ containing them as
rectilinear generators. As $p^i$ tend to $p^0$, the quadric ${\bf Q}$ tends to
a
limiting quadric, the so-called Lie quadric of the surface $M^2$ at the point
$p^0$. Even though this construction depends on the initial choice of either
the $x$-
or the $y$-asymptotic line through $p^0$, the resulting quadric ${\bf Q}$ is
independent
of that choice. Thus, we arrive at a two-parameter family of quadrics
associated with the surface $M^2$. In terms of the Wilczynski tetrahedral, the
parametric equation for ${\bf Q}$ is of the form \cite[p.311]{Bol1}
$$
{\bf Q}={\mbox{\boldmath $\eta$}}+\mu {\bf r}_1+\nu {\bf r}_2+\mu\nu {\bf r},
$$
where $\mu,\nu$ are parameters.

Now, in the neighbourhood of a generic point $p^0$ on $M^2$, the envelopes of
the family of Lie quadrics consist of the surface $M^2$ itself and four, in
general, distinct sheets.
The case of projectively minimal surfaces is characterized by the additional
requirement that the asymptotic lines on all these sheets correspond to the
asymptotic lines of the surface $M^2$ itself \cite{Thomsen2}.
Moreover, for projectively minimal
surface all four sheets of the envelope will be
projectively minimal as well \cite{Mayer}.
In a sense, it is natural to
call the family of Lie quadrics with this property a W-congruence of quadrics.
Surfaces of Godeaux-Rozet are characterized by the
degenerate case of two distinct sheets while Demoulin surfaces are present
if all four sheets coincide.
Surfaces of Godeaux-Rozet and Demoulin
have been investigated extensively in~\cite{Demoulin,Godeaux,Rozet},
see also \cite{Fer8}.

\section{Congruences W}

There exists an important class of transformations in projective differential
geometry which leave the system (\ref{r}) form-invariant. These are
transformations
generated by congruences W. Here we briefly recall this construction following
\cite{Jonas1}, \cite{Finikov37}, \cite{Eisenhart}.

Let $M^2$ be a surface with the radius-vector ${\bf r}$ satisfying (\ref{r}).
Let the functions ${\cal U}$ and ${\cal V}$ satisfy the Dirac equation
\begin{equation}
\begin{array}{c}
{\cal U}_x=\beta \ {\cal V} \\
{\cal V}_y=\gamma \ {\cal U}
\end{array}
\label{W}
\end{equation}
where $\beta$ and $\gamma$ are the same as in (\ref{r}). Let us consider a
surface
$\tilde M^2$ with the radius-vector
${\bf r}'$ given by the formula
\begin{equation}
{\bf r}'={\cal V} \ {\bf r}_1- {\cal U} \ {\bf r}_2+
\frac{1}{2}\left( {\cal V}
 \ \frac{\gamma_x}{\gamma} -
{\cal U} \ \frac{\beta_y}{\beta}-{\cal V}_x+{\cal U}_y\right) ~ {\bf r}
\label{newtilder}
\end{equation}
In order to write down the equations for ${\bf r}'$ it is convenient to
introduce certain quantities which are combinations of
${\cal U}, {\cal V}$ and their derivatives. First of all, we define
${\cal A}$ and  ${\cal B}$ by the formulae
$$
{\cal U}_y=\frac{\beta_y}{\beta} \ {\cal U}+{\cal A}, ~~~~
{\cal V}_x=\frac{\gamma_x}{\gamma} \ {\cal V}+{\cal B},
$$
(in fact, we are copying equations (\ref{UAPVBQ}) for the Pl\"ucker
coordinates).
The compatibility conditions ${\cal U}_{xy}={\cal U}_{yx}$ and
${\cal V}_{xy}={\cal V}_{yx}$ imply
$$
{\cal A}_x = k \ {\cal U}, ~~~~ {\cal B}_y=l \ {\cal V},
$$
where $l$ and $k$ are the same as in (\ref{klab}). Let us introduce
${\cal P}$ and ${\cal Q}$ by the formulae
$$
{\cal A}_y=a \ {\cal U}+{\cal P}, ~~~~ {\cal B}_x=b \ {\cal V}+{\cal Q}.
$$
Then the compatibility conditions imply
$$
{\cal P}_x=-\beta a \ {\cal V}+k \ {\cal A}, ~~~~
{\cal Q}_y=-\gamma b \ {\cal U}+l \ {\cal B}.
$$
Finally, we introduce the quantities $H$ and $K$ via
$$
{\cal P}_y=a \ {\cal A}-\frac{\beta_y}{\beta} \ {\cal P} -\gamma b \ {\cal V}+
\gamma \ {\cal Q}+K, ~~~~
{\cal Q}_x=b \ {\cal B}-\frac{\gamma_x}{\gamma} \ {\cal Q} -\beta a \ {\cal U}+
\beta \ {\cal P}+H,
$$
so that the compatibility conditions imply
$$
H_y=-\beta \ K, ~~~~ K_x=-\gamma \ H.
$$
Equations for ${\cal U}, {\cal A}, {\cal P}, {\cal V}, {\cal B}, {\cal Q}, H,
K$
can be rewritten in the matrix form
\begin{equation}
\begin{array}{c}
\left(\begin{array}{c}
{\cal U}\\
{\cal A}\\
{\cal P}\\
{\cal V}\\
{\cal B}\\
{\cal Q}\\
H\\
K
\end{array}\right)_x=
\left(\begin{array}{cccccccc}
0 & 0 & 0 & \beta & 0 & 0 & 0 & 0\\
k & 0 & 0 & 0 & 0 & 0 & 0 & 0\\
0 & k & 0 & -\beta a & 0 & 0 & 0 & 0\\
0 & 0 & 0 & \frac{\gamma_x}{\gamma} & 1 & 0 & 0 & 0\\
0 & 0 & 0 & b & 0 & 1 & 0 & 0\\
-\beta a & 0 & \beta & 0 & b &-\frac{\gamma_x}{\gamma} & 1 & 0\\
{ * } & * & * & * & * & * & * & * \\
0 & 0 & 0 & 0 & 0 & 0 & -\gamma & 0
\end{array}\right)
\left(\begin{array}{c}
{\cal U}\\
{\cal A}\\
{\cal P}\\
{\cal V}\\
{\cal B}\\
{\cal Q}\\
H\\
K
\end{array}\right)\\
\ \\
\left(\begin{array}{c}
{\cal U}\\
{\cal A}\\
{\cal P}\\
{\cal V}\\
{\cal B}\\
{\cal Q}\\
H\\
K
\end{array}\right)_y=
\left(\begin{array}{cccccccc}
\frac{\beta_y}{\beta} & 1 & 0 & 0 & 0 & 0 & 0 & 0\\
a & 0 & 1 & 0 & 0 & 0 & 0 & 0\\
0 & a & -\frac{\beta_y}{\beta} & -\gamma b & 0 & \gamma & 0 & 1\\
\gamma & 0 & 0 & 0 & 0 & 0 & 0 & 0\\
0 & 0 & 0 & l & 0 & 0 & 0 & 0\\
-\gamma b & 0 & 0 & 0 & l & 0 & 0 & 0\\
 0 & 0 & 0 & 0 & 0 & 0 & 0 & -\beta\\
{ * } & * & * & * & * & * & * & *
\end{array}\right)
\left(\begin{array}{c}
{\cal U}\\
{\cal A}\\
{\cal P}\\
{\cal V}\\
{\cal B}\\
{\cal Q}\\
H\\
K
\end{array}\right)
\end{array}
\label{W1}
\end{equation}
where the elements $*$ are not yet specified. Equations (\ref{W1}) reduce to
(\ref{UAPVBQ}) under the reduction $H=K=0$.
In what follows we will also need the quantity
$$
S={\cal Q} \ {\cal V}-{\cal P} \ {\cal U}+\frac{{\cal A}^2-{\cal B}^2}{2},
$$
which, in view of (\ref{W1}), satisfies the equations
$$
S_x=H \ {\cal V}, ~~~~ S_y=-K \ {\cal U}.
$$
Now a direct calculation gives:
\begin{equation}
{\bf r}=-2 \frac{{\cal V}}{S} \  {\bf r}'_x-
2 \frac{{\cal U}}{S}\ {\bf r}'_y+\frac{1}{S}
({\cal A}+{\cal B}+\frac{\gamma_x}{\gamma}{\cal V}+\frac{\beta_y}{\beta}{\cal
U})
\ {\bf r}'.
\label{oldr}
\end{equation}
Equations (\ref{newtilder}) and (\ref{oldr}) immediately imply that the line
${\bf r}\wedge  {\bf r}'$ joining the corresponding points
${\bf r}$ and ${\bf r}'$ is tangent to both surfaces $M^2$ amd $\tilde M^2$,
which are thus the focal surfaces of the line congruence
${\bf r}\wedge  {\bf r}'$. Moreover, the formulae
\begin{equation}
\begin{array}{c}
{\bf r}'_{xx}=\frac{S_x}{S}\  {\bf r}'_x+
\left (\frac{S_x}{S} \frac{{\cal U}}{{\cal V}}-\beta \right)\ {\bf r}'_y+
\frac{1}{2}\left( V+\beta_y-\frac{S_x}{S{\cal V}}
({\cal A}+{\cal B}+\frac{\gamma_x}{\gamma}{\cal V}+\frac{\beta_y}{\beta}{\cal
U})
\right) {\bf r}', \\
\ \\
{\bf r}'_{yy}=\frac{S_y}{S}\  {\bf r}'_y+
\left (\frac{S_y}{S} \frac{{\cal V}}{{\cal U}}-\gamma \right)\  {\bf r}'_x+
\frac{1}{2}\left( W+\gamma_x-\frac{S_y}{S{\cal U}}
({\cal A}+{\cal B}+\frac{\gamma_x}{\gamma}{\cal V}+\frac{\beta_y}{\beta}{\cal
U})
\right)  {\bf r}',
\end{array}
\label{W2}
\end{equation}
(which are the result of quite a long calculation)
demonstrate that $x, y$ are asymptotic coordinates on the transformed surface
$\tilde M^2$ as well, so that the congruence ${\bf r}\wedge {\bf r}'$
preserves the asymptotic parametrization of it's focal surfaces. Such
congruences play
a central role in projective differential geometry and are known as
the congruences W. By a construction, a congruence W with one given focal
surface
$M^2$ is uniquely determined by a solution ${\cal U}, {\cal V}$ of the linear
Dirac equation (\ref{W}). Normalising the vector ${\bf r}'$ as follows:
${\bf r}'={\sqrt S}~\tilde {\bf r}$ we can rewrite equations (\ref{W2}) in the
canonical form (\ref{r}) in terms of $\tilde {\bf r}$:
\begin{equation}
\begin{array}{c}
\tilde {\bf r}_{xx}=\tilde \beta \ \tilde {\bf r}_y+
\frac{1}{2}(\tilde V-\tilde \beta_y) \ \tilde {\bf r} \\
\ \\
\tilde {\bf r}_{yy}=\tilde \gamma \ \tilde {\bf r}_x+
\frac{1}{2}(\tilde W-\tilde \gamma_x) \ \tilde {\bf r}
\end{array}
\label{tilder}
\end{equation}
where the transformed coefficients $\tilde \beta, \ \tilde \gamma, \ \tilde V,
\
\tilde W$ are given by the formulae
\begin{equation}
\begin{array}{c}
\tilde \beta=\frac{S_x}{S} \frac{{\cal U}}{{\cal V}}-\beta=
 \frac{H {\cal U}}{S}-\beta, ~~~~~~
\tilde \gamma=\frac{S_y}{S} \frac{{\cal V}}{{\cal U}}-\gamma=
- \frac{K {\cal V}}{S}-\gamma,
 \\
\ \\
\tilde V = V-\frac{S_x}{S} \frac{{\cal V}_x}{{\cal V}}
+\frac{3}{2}(\frac{S_x}{S})^2-\frac{S_{xx}}{S}, ~~~~~
\tilde W = W-\frac{S_y}{S} \frac{{\cal U}_y}{{\cal U}}
+\frac{3}{2}(\frac{S_y}{S})^2-\frac{S_{yy}}{S},
\end{array}
\label{W3}
\end{equation}
(we point out the usefull identity
$\tilde \beta \tilde \gamma = \beta \gamma - (\ln S)_{xy}$).
We will call the surface $\tilde M^2$ a W-transform of the surface $M^2$.
Congruences W provide a standard tool for constructing
B\"acklund transformations. Suppose we are given a
class of surfaces specified by a certain  constraint
between $\beta, \ \gamma, \ V, \ W$.
Let us try to find a congruence W such that
the second focal surface will also belong to the same class.
This requirement imposes
additional restrictions on the functions ${\cal U}$ and ${\cal V}$,
which usually turn to be linear  and, moreover, contain an arbitrary
constant parameter, so that equations (\ref{W1})
become a "Lax pair" for the class of surfaces under study.
Since the Dirac equation (\ref{W}) is a part of this Lax pair, it is not
surprising
that surfaces in projective geometry are closely related to the DS hierarchy.
Particularly interesting classes of surfaces correspond to
reductions of the Dirac operator which are quite familiar from the modern
soliton theory. These are

isothermally-asymptotic surfaces ($\beta=\gamma$);

surfaces $R_0$ ($\beta=1$ or $\gamma=1$),

surfaces R ($\beta_y=\gamma_x$);

surfaces of Jonas ($\beta_x=\gamma_y$), etc.

\bigskip

{\bf Example 1. B\"acklund transformation for isothermally asymptotic
surfaces}.
Let us require that both surfaces $M^2$ and $\tilde M^2$ are isothermally
asymptotic,
that is, $\beta=\gamma$ and $\tilde \beta = \tilde \gamma$.
For that purpose it is sufficient to choose
$$
H=\lambda \ {\cal V}, ~~~~ K=-\lambda \ {\cal U}, ~~~~ \lambda = const,
$$
which, upon the substitution in (\ref{W1}), results in the Lax pair
(\ref{UAPVBQ1}). The transformation
$$
\tilde \beta=\lambda \frac{{\cal U}{\cal V}}{S}-\beta
$$
provides thus a B\"acklund transformation of the stationary mVN equation.

\bigskip

{\bf Example 2. B\"acklund transformation for surfaces $R_0$}.
Here we require $\beta=\tilde \beta=1$, implying $H{\cal U}=2S$.
Differentiation by
$x$ and $y$ produces further constraints $H_x{\cal U}=H{\cal V}$ and
$K{\cal U}+H{\cal A}=0$, respectively. These restrictions can be identically
satisfied
if we assume
$$
H=-\lambda \ {\cal U}, ~~~~ K=\lambda \ {\cal A},
$$
 and impose the compatible quadratic constraint
$$
\lambda \ {\cal U}^2+2 \ S=0.
$$

\bigskip

{\bf Example 3. B\"acklund transformation for surfaces $R$}.
Here we require $\beta_y=\gamma_x$ and
$\tilde \beta_y=\tilde \gamma_x$ implying $(H{\cal U}/S)_y+(K{\cal V}/S)_x=0$.
To satisfy this condition, it is sufficient to choose
$$
H=\lambda \left(\frac{\gamma_x}{\gamma} \ {\cal V}+{\cal B}-
\beta \ {\cal U}\right), ~~~~
K=\lambda \left(\frac{\beta_y}{\beta} \ {\cal U}+{\cal A}-
\gamma \ {\cal V}\right)
$$
and to impose the compatible quadratic constraint
$$
\lambda \ \left({\cal U}^2-{\cal V}^2\right)+2 \ S=0.
$$
In a different gauge this B\"acklund transformation has been set down by Jonas
in
\cite{Jonas1}.

\bigskip

{\bf Example 4. B\"acklund transformation for surfaces of Jonas}.
Here we require $\beta_x=\gamma_y$ and
$\tilde \beta_x=\tilde \gamma_y$ implying $(H{\cal U}/S)_x+(K{\cal V}/S)_y=0$.
This restiction is satisfied if $H$ and $K$ satisfy  the equations
$$
H_x= \lambda \ {\cal V} -\gamma \ K, ~~~~
K_y= \lambda \ {\cal U} -\beta  \ H,
$$
along with the additional quadratic constraint
$$
2\lambda  \ S+K^2-H^2=0.
$$
With these $H, K$ equations (\ref{W1}) transform into the $8\times 8$ linear
problem (\ref{8x8}) which is compatible with the above quadratic constraint.
This B\"acklund transformation has
been set down by Jonas in \cite{Jonas2}.

\bigskip

{\bf Example 5. B\"acklund transformation for surfaces with one family of
asymptotic lines in linear complexes}. This class of surfaces is specified by
the condition $k=0$  in equations (\ref{GC2}), implying
$(\ln \beta) _{xy}=\beta \gamma$.
Geometrically, this  means that the Pl\"ucker image ${\cal U}=
{\bf r}\wedge {\bf r}_1$ of any $x$-asymptotic curve of the surface $M^2$
is a hyperplane curve in $P^5$ (one can similarly require $l=0$ or
$(\ln \gamma) _{xy}=\beta \gamma$ implying the same property for
 y-asymptotic curves). Indeed, the condition $k=0$ implies
${\cal A}_x=0$ (see (\ref{UAPVBQ})), so that the vector ${\cal A}$
is constant along any $x$-asymptotic line. The conditions
$({\cal A}, {\cal U})=({\cal A}, {\cal U}_x)=...=0$ (see (\ref{table})) means
that
the curve ${\cal U}$ lies in a hyperplane in $P^5$
which is orthogonal to ${\cal A}$.
We recall, that by a definition a linear complex is a 3-parameter family of
straight lines in $P^3$, which corresponds to a hyperplane in $P^5$ under the
Pl\"ucker embedding. Surfaces with asymptotic lines in linear complexes have
been
extensively studied in projective differential geometry, see e.g.
\cite{Jonas3},
\cite{Keraval}, \cite{Finikov50}
\cite{Fubini} and references therein. The most important
geometric property of these surfaces is the existence of a W-congruence,
mapping a surface $M^2$ onto a ruled surface $\tilde M^2$ for which
$\tilde \beta =0$. In this case the $x$-asymptotic curves of the transformed
surface are just straight lines. This construction rectifies all $x$-asymptotic
curves of the surface
$M^2$ simultaneosly and in fact linearises equations (\ref{GC2}) with the
constraint $k=0$. Requiring in (\ref{W3}) $\tilde \beta=0$, or, equivalently,
$H {\cal U}=\beta S$, and differentiating this constraint with respect
to $x$ and $y$, we arrive at the additional constraints
$$
H_x=\frac{\beta_x}{\beta} H ~~~~~ {\rm and} ~~~~~ {\cal A}=0,
$$
respectively. Taking into account (\ref{W1}), the last condition implies
$$
{\cal P}=-a \ {\cal U}, ~~~~~ K=-(a_y+2 a \frac{\beta_y}{\beta}) \ {\cal U}+
\gamma b \ {\cal V} -\gamma \ {\cal Q},
$$
so that equations (\ref{W1}) reduce to
$$
\begin{array}{c}
\left(\begin{array}{c}
{\cal U}\\
{\cal V}\\
{\cal B}\\
{\cal Q}\\
H\\
\end{array}\right)_x=
\left(\begin{array}{ccccc}
0 &  \beta  & 0 & 0 & 0\\
0 & \frac{\gamma_x}{\gamma} & 1  & 0 & 0\\
0 &  b & 0 & 1 & 0 \\
-2 \beta a & 0  & b &-\frac{\gamma_x}{\gamma} & 1 \\
0 & 0 & 0 & 0 & \frac{\beta_x}{\beta}
\end{array}\right)
\left(\begin{array}{c}
{\cal U}\\
{\cal V}\\
{\cal B}\\
{\cal Q}\\
H\\
\end{array}\right)\\
\ \\
\left(\begin{array}{c}
{\cal U}\\
{\cal V}\\
{\cal B}\\
{\cal Q}\\
H\\
\end{array}\right)_y=
\left(\begin{array}{ccccc}
\frac{\beta_y}{\beta}  & 0 & 0 & 0 & 0\\
\gamma & 0 & 0 & 0 & 0\\
0 &  l  & 0 & 0 & 0\\
-\gamma b & 0 & l  & 0 & 0\\
\beta a_y+2a\beta_y & -\beta \gamma b & 0 &\beta \gamma & 0
\end{array}\right)
\left(\begin{array}{c}
{\cal U}\\
{\cal V}\\
{\cal B}\\
{\cal Q}\\
H\\
\end{array}\right)
\end{array}
$$
This system is compatible and preserves the constraint $H {\cal U} = \beta S$.
Thus, it defines a 3-parameter family of transformations W mapping a surface
$M^2$ with $x$-asymptotic lines in linear complexes onto a ruled surface
$\tilde M^2$. Conversely, applying an arbitrary transformation W to a ruled
surface, we obtain a surface with one family of asymptotic lines in
linear complexes (those corresponding to  rectilinear generators of
a ruled surface).

\bigskip

{\bf Example 5. B\"acklund transformation for surfaces with both families of
asymptotic lines in linear complexes}. This class of surfaces is specified by
the condition $k=l=0$  implying a pair of the Liouville equations for
$\beta, \gamma$:
$$
(\ln \beta) _{xy}=\beta \gamma, ~~~~~ (\ln \gamma) _{xy}=\beta \gamma.
$$
In this case one can find a transformation W, which maps  $M^2$
onto a surface $\tilde M^2$ for which $\tilde \beta = \tilde \gamma =0$, that is,
onto a quadric. Thus, it simultaneously rectifies all asymptotic curves.
In order to construct such a W-transformation we have to take
$$
{\cal A}={\cal B}=0, ~~~~ {\cal P}=-a \ {\cal U}, ~~~~ {\cal Q}=-b \ {\cal V},
$$
$$
H=-(b_x+2b\frac{\gamma_x}{\gamma}) \ {\cal V} +2\beta a \ {\cal U}, ~~~~~
K=-(a_y+2a\frac{\beta_y}{\beta}) \ {\cal U} +2\gamma b \ {\cal V},
$$
implying the following equations for ${\cal U}, {\cal V}$
$$
{\cal U}_x=\beta \ {\cal V}, ~~~~ {\cal U}_y=\frac{\beta_y}{\beta} \ {\cal U},
~~~~
{\cal V}_x=\frac{\gamma_x}{\gamma} \ {\cal V}, ~~~~
{\cal V}_y=\gamma \ {\cal U}.
$$
together with the compatible quadratic constraint
$$
\beta \gamma \ (a \ {\cal U}^2+ b \ {\cal V}^2) = (\gamma b_x + 2 b \gamma _x)
\ {\cal U}{\cal V}.
$$
These ${\cal U}, {\cal V}$ define a unique transformation W, mapping our
surface to a quadric.
Conversely, applying an arbitrary W-transformation to a quadric, we obtain
a surface with both families of asymptotic lines in linear complexes.

\section{Acknowledgements}

I would like to thank A.I.~Bobenko, K.R.~Khusnutdinova, B.G.~Konopelchenko,
U.~Pinkall and W.K.~Schief for their interest and useful remarks.

This research was partially supported by INTAS 96-0770 and
the Alexander von Humboldt Foundation.

\end{document}